\documentclass[10pt]{article}
\font\bfone=cmbx12 scaled \magstep1

\begin{document}

\centerline{\bfone Lattice points in large regions and related
arithmetic functions:}\smallskip\centerline{\bfone Recent
developments in a very classic topic}\bigskip

\centerline{\bf A.~Ivi\'c (Belgrade),}\smallskip \centerline{\bf
E.~Kr\"atzel, M.~K\"uhleitner, and W.G.~Nowak (Vienna)}\bigskip

\centerline{\it Dedicated to Professor Wolfgang Schwarz on the
occasion of his 70th birthday}\bigskip

\def\ssk{\par\smallskip}
\def\msk{\par\medskip}
\def\bsk{\par\bigskip}
\def\cen{\centerline}
\def\Int{\int\limits}

\def\dint #1 {
\quad  \setbox0=\hbox{$\displaystyle\int\!\!\!\int$}
  \setbox1=\hbox{$\!\!\!_{#1}$}
  \vtop{\hsize=\wd1\centerline{\copy0}\copy1} \quad}

\font\bfone=cmbx10 scaled\magstep1 %
\font\boldmas=msbm10                  
\def\Bbb#1{\hbox{\boldmas #1}}        
\def\Z{{\Bbb Z}}                        
\def\N{{\Bbb N}}                        
\def\Q{{\Bbb Q}}
\def\R{{\Bbb R}}
\def\C{{\Bbb C}}
\font\boldmasi=msbm10 scaled 700      
\def\Bbbi#1{\hbox{\boldmasi #1}}      
\font\boldmas=msbm10                  
\def\Bbb#1{\hbox{\boldmas #1}}        
\def\Zi{{\Bbbi Z}}                      
\def\Ci{{\Bbbi C}}                      
\def\Ni{{\Bbbi N}}                      
\def\Pi{{\Bbbi P}}                      
\def\Qi{{\Bbbi Q}}                      
\def\Ri{{\Bbbi R}}

The branch of analytic number theory which is concerned with the
number of integer points in large domains, in the Euclidean plane
as well as in space of dimensions $\ge3$, has a long and very
prolific history which reaches back to classic works of E.~Landau,
J.G.~Van der Corput, G.~Vorono\"{\i}, G.H.~Hardy, to mention just
a few. Enlightening accounts of the developments of this theory
can be found in the monographs of F.~Fricker \cite{fricker},
E.~Kr\"atzel \cite{kr-lp},\cite{kr-anafu}, and M.~Huxley
\cite{Hu1}. Although the major problems in this field are quite
old (and yet unsolved), there have been a lot of amazing new
developments in recent times. This survey article attempts to give
an overview of the state-of-art in this field, with an emphasis on
results established during the last few years. Furthermore, we are
aiming to point out the intrinsic relationships between lattice
point quantities and certain arithmetic functions which arose
originally from number theory originally without any geometric
motivation.

\subsubsection*{1. Lattice points in a large circle. The problem of Gau\ss.}
Starting from the arithmetic function viewpoint, let

$$r(n) := \#\{(m,k)\in\Z^2:\ m^2+k^2 = n \ \}$$
denote the number of ways to represent the integer $n\ge0$ as a
sum of two squares of integers. While the maximal size of $r(n)$
is estimated, e.g., in E. Kr\"atzel \cite{kr}, and an explicit
formula may be found in E.~Hlawka, J.~Schoi\ss engeier, and
R.~Taschner \cite{hst}, p.~60, its average order is described by
the Dirichlet summatory function
$$A(x) = \sum_{0\le n\le x} r(n)\,,$$
where $x$ is a
large real variable. In geometric terms, $A(x)$ is the number of
lattice points $(m,k)\in\Z^2$ in a compact, origin-centered
circular disc ${\bf C}(\sqrt{x})$ of radius $\sqrt{x}$. Obviously,
$A(x)$ equals asymptotically the area of ${\bf C}(\sqrt{x})$,
i.e., $\pi x$. Already C.F.~Gau\ss\ (1777 - 1855) made the
following more precise but still extremely simple observation: If
we associate to each lattice point $(m,k)\in\Z^2$ the square
$$Q(m,k) = \{ (u,v)\in\R^2:\ |u-m|\le{\textstyle{1\over2}},\
|v-k|\le{\textstyle{1\over2}}\ \}\,,$$ then obviously
$${\bf C}\left(\sqrt{x}-{\sqrt{2}\over2}\right) \subset
\bigcup_{(m,k)\in{\bf C}(\sqrt{x})}Q(m,k)\subset{\bf
C}\left(\sqrt{x}+{\sqrt{2}\over2}\right)\,.$$ Evaluating the
areas, we see that, for every $x>{1\over2}$, $$
\pi\left(\sqrt{x}-{\sqrt{2}\over2}\right)^2 < A(x) <
\pi\left(\sqrt{x}+{\sqrt{2}\over2}\right)^2\,,$$ hence the {\it
lattice point discrepancy } of ${\bf C}(\sqrt{x})$ (traditionally
often mis-translated as "lattice rest") $$P(x):= A(x)-\pi x
\eqno{(1.1)}$$ satisfies $$ |P(x)| < \pi\sqrt{2x}+{\pi\over2}
\eqno{(1.2)}$$ which is usually simply written as
$P(x)=O(\sqrt{x})$ or $P(x)\ll\sqrt{x}$. \msk Before reporting on
stronger results, we shall try to shed some light on the
mathematical background of all more sophisticated approaches to
the problem. The na\"{\i}ve idea is to apply Poisson's formula to
evaluate the sum
$$A(x) = \sum_{(m,k)\in\Zi^2}{\bf c}_x(m,k)\,,$$
where ${\bf c}_x$ denotes the indicator function of ${\bf
C}(\sqrt{x})$. To get rid of the discontinuity at the boundary of
${\bf C}(\sqrt{x})$, which would create subtle questions of
convergence, it is convenient to consider the integrated form
$$ \int\limits_0^x A(t)\, {\rm d}t = \int\limits_0^x \left(\sum_{m^2+k^2\le t} 1\right) {\rm d}t = \sum_{m^2+k^2\le
x}(x-m^2-k^2) = $$ $$ = \sum_{(m,k)\in\Zi^2}\ \dint{u^2+v^2\le x}
(x-u^2-v^2)e^{2\pi i(mu+kv)}\, {\rm d}(u,v)\,. $$ Since the double
integral can be evaluated explicitly by means of the Bessel
function $J_2$, this gives
$$  \int\limits_0^x A(t)\, {\rm d}t = {\pi\over2}\,x^2 + {x\over\pi}
\sum_{n=1}^\infty {r(n)\over n}\,J_2\left(2\pi\sqrt{nx}\right)\,.
\eqno{(1.3)}  $$ With some care, one can verify that this identity
may be differentiated term by term (cf., e.g., E. Kr\"atzel
\cite{kr-lp}, p.~126, Th.~3.12), yielding (for $x$ not an integer)
$$ P(x) = \sqrt{x} \sum_{n=1}^\infty r(n)n^{-1/2}
J_1(2\pi\sqrt{nx})\,, \eqno{(1.4)}  $$ which is called {\it
Hardy's identity}. Since this latter series is only conditionally
convergent, it is convenient for applications to have  a sharp
truncated version at hand. A result of this kind can be found as
Lemma 1 in A. Ivi\'c \cite{I8}: For $x\ge1$, $x\notin\Z$, and
$x\le M\le x^A$, $A>1$ some fixed constant,
$$P(x) = \sqrt{x} \sum_{1\le n< M} r(n)n^{-1/2}
J_1(2\pi\sqrt{nx}) +\qquad\qquad\qquad\qquad\qquad\qquad\quad$$
$$+ O\left(\min\left(x^{5/4}M^{-1/2}+x^{1/2+\epsilon}M^{-1/2}\Vert
x\Vert^{-1} + x^{1/4}M^{-1/4},\ x^\epsilon\right)\right)\,,
\eqno{(1.5)}$$ where $\Vert\cdot\Vert$ denotes the distance from
the nearest integer, and $\epsilon>0$ is arbitrary. Moreover,
using the well-known Bessel function asymptotics
$$ J_1(2\pi\sqrt{nx}) = {1\over\pi}(nx)^{-1/4} \cos(2\pi\sqrt{nx}-{\textstyle{3\over4}}\pi) +
O\left((nx)^{-3/4}\right)\,, $$ one arrives at a
representation
$$ P(x) = {1\over\pi}\,x^{1/4}\sum_{1\le n<M} {r(n)\over n^{3/4}}\,
\cos(2\pi\sqrt{nx}-{\textstyle{3\over4}}\pi) \ + \hbox{remainder
terms}\,.\eqno{(1.6)} $$ In a sense, all of the deep
investigations on the order of $P(x)$ which have been done during
the past century, are based on results like (1.6),
though often somewhat shrouded by technicalities.

\paragraph*{1.1.~O-estimates in the
circle problem. } For instance, to derive upper bounds for this
lattice point discrepancy, one is led to consider exponential sums
$$ \sum_{m,k}(m^2+k^2)^{-3/4}\,e^{2\pi i\sqrt{x}\sqrt{m^2+k^2}}\,.
$$ The twentieth century saw the development of increasingly
complicated methods to deal with such sums. A historical survey of
the results obtained can be found in E. Kr\"atzel \cite{kr-lp}.
They started with W.~Sierpi\'nski's \cite{si} bound $$ P(x) =
O\left(x^{1/3}\right)
$$ and in the 1980's culminated in the work of G.A.~Kolesnik \cite{ko},
who had
$$ P(x) = O\left(x^{139/429 +\epsilon}\right)  $$ as his sharpest result.
\ssk During the last one-and-a-half decades, M.~Huxley, starting
from the ideas due to E. Bombieri -- H. Iwaniec \cite{BI} and
ideas due to H.~Iwaniec and C.J.~Mozzochi \cite{IM}, devised a
substantially new approach which he called the {\it Discrete
Hardy-Littlewood method}. His strongest result was published in
2003 \cite{Hu2} and reads
$$  P(x) = O\left(x^{131/416} (\log x)^{18637/8320}\right)\,,  \eqno{(1.7)}
$$ thereby improving upon his 1993 \cite{hu} bound $O(x^{23/73+\epsilon})$.
Note that ${131\over416}=0.3149038\dots$,
${23\over73}=0.315068\dots$, while Kolesnik's
${139\over429}=0.324009\dots$.

\paragraph*{\bf 1.2.~Lower bounds for the
lattice point discrepancy of a circle. } As (1.6) suggests, a
stronger estimate than $P(x)=O(x^{1/4})$ cannot be true. In fact,
G.H.~Hardy \cite{ha}, \cite{Ha1} proved that\footnote{Let
us recall briefly the different $\Omega$-symbols: For real
functions $F$ and $G>0$, we mean by \hbox{$F(x)=\Omega_+(G(x))$}
that $\limsup(F(x)/G(x))>0$, as $x\to\infty$. Similarly,
$F(x)=\Omega_-(G(x))$ stands for $-F(x)=\Omega_+(G(x))$. Further,
$F(x)=\Omega_\pm(G(x))$ says that both of these assertions are
true, and $F(x)=\Omega(G(x))$ means that at least one of them
applies, i.e., $F(x)\ne o(G(x))$.}
$$P(x)=\Omega_\pm(x^{1/4})\,,\qquad P(x)=\Omega_-(x^{1/4}(\log
x)^{1/4})\,, \eqno(1.8)  $$  adding the unsubstantiated claim that
$P(x)=\Omega_\pm(x^{1/4}(\log x)^{1/4})$. The second part of (1.8)
was not improved until more than six decades later by J.L.~Hafner
\cite{Hf1}, who obtained\footnote{We write throughout $\log_j$ for
the $j$-fold iterated natural logarithm. Thus $\log_2=\log\log$,
and so on.}
$$ P(x)=\Omega_-\left(x^{1/4}(\log x)^{1/4} (\log_2x)^{(\log2)/4}
\exp(-c(\log_3x)^{1/2})\right)\,,\qquad c>0\,. \eqno{(1.9)}  $$
All estimates like (1.8) and (1.9) start from (a variant of)
formula (1.6) and are based on the idea of finding an unbounded
sequence of $x$-values for which the terms
$\cos(2\pi\sqrt{nx}-{\textstyle{3\over4}}\pi)$ "essentially all
pull into the same direction". That is, the numbers $\sqrt{nx}$
should be close to an integer, which is achieved by Dirichlet's
theorem on Diophantine approximation (see e.g., Lemma 9.1 of
\cite{I13}). In Hafner's work, this task is done only for those
$n$ for which $r(n)$ is comparatively large; this idea slightly
improves the estimate. \ssk As for $\Omega_+$-results,
K.~Corr\'adi and I.~K\'atai \cite{CK} proved that
$$ P(x) = \Omega_+\left(x^{1/4}\exp\left(c'(\log_2x)^{1/4}
(\log_3x)^{-3/4}\right)\right)\qquad (c'>0)\,,  \eqno{(1.10)} $$
thereby refining previous work of A.E.~Ingham \cite{ing} and of
K.S.~Gangadharan \cite{Ga}. (Note that the factor of $x^{1/4}$
here is less than any power of $\log x$. This inherent asymmetry
arises from the fact that the deduction of (1.10) is based on sort
of a quantitative version of Kronecker's approximation theorem
which is much weaker than Dirichlet's.) \ssk The strongest
$\Omega$-bound for $P(x)$ available to date was established  in
2003 by K.~Soundararajan \cite{So}. His ingenious new approach
allows us to restrict the application of
Dirichlet's theorem to a still smaller set of integers $n$. His
method can be summarized in the following fairly general
statement. \msk\noindent{\bf Soundararajan's Lemma {\rm
\cite{So}}.}\quad{\it Let $(f(n))_{n=1}^\infty$ and
$(\lambda_n)_{n=1}^\infty$ be sequences of non\-negative real
numbers, $(\lambda_n)_{n=1}^\infty$ non-decreasing, and
$\sum_{n=1}^\infty f(n)<\infty$. Let $L\ge2$ and $Y$ be integer
parameters, $\gamma\in\R$ fixed. Put
$$ S(t) = \sum_{n=1}^\infty f(n)\cos(2\pi\lambda_n t + \gamma)\,. $$
Suppose further that ${\cal M}$ is a finite set of positive
integers, such that $\{\lambda_m:\ m\in{\cal M}\ \}
\subset[{1\over2}\lambda_Y,{3\over2}\lambda_Y]$. Then, for any
real $T\ge2$, there exists some $t\in[{1\over2} T, (6L)^{|{\cal
M}|+1}\,T]$ such that
$$ \left|S(t)\right|\ \ge\ {1\over8}\sum_{m\in{\cal M}}f(m) -
{1\over L-1}\sum_{n:\ \lambda_n\le2\lambda_Y}f(n) -
{4\over\pi^2T\lambda_Y}\sum_{n=1}^\infty f(n)\,. $$ Furthermore,
if $\gamma=w\pi$, $w\in\Z$, then on the left hand side $|S(t)|$
can be replaced by $(-1)^w\,S(t)$.} \bsk When applied to (1.6),
this lemma yields
$$ P(x) = \Omega\left(x^{1/4}(\log x)^{1/4}(\log_2 x)^{(3/4)(2^{1/3}-1)}
(\log_3 x)^{-5/8}\right)\,. \eqno{(1.11)} $$ To assess the
refinement in the exponent of the $\log_2 x$-factor, note that
${1\over4}\log2=0.1732\dots$, while
${3\over4}(2^{1/3}-1)=0.1949\dots$. We remark that Soundararajan's
method does not allow us to replace the $\Omega$-symbol in this
assertion by $\Omega_+$ or $\Omega_-$.

\paragraph*{1.3.~The mean square in the circle problem. } The results reported so far already
suggest that
$$ \inf\{\lambda:\ P(x) \ll_\lambda\ x^\lambda\ \} = {\textstyle{1\over4}}\,.
\eqno{(1.12)}  $$ To prove (or disprove) this conjecture is
usually called the Gaussian circle problem in the strict sense. In
favor of this hypothesis, there are quite precise mean-square
asymptotic formulas of the shape
$$ \Int_0^X (P(x))^2\,  {\rm d}x = C\,X^{3/2} + Q(X)\,, \eqno{(1.13)}  $$
with
$$  C = {1\over3\pi^2}\sum_{n=1}^\infty r^2(n) n^{-3/2} =
{16\over3\pi^2}\,
{\zeta_{\Qi(i)}^2({3\over2})\over\zeta(3)}(1+2^{-3/2})^{-1}
 \approx 1.69396\,. \eqno{(1.14)}  $$
The estimation of this new remainder term $Q(X)$ has been subject
of intensive research, using more and more ingenuity. We mention
the results of H.~Cram\'er \cite{cr}: $Q(X)=O( X^{5/4+\epsilon})$,
E.~Landau \cite{la3}: $Q(X)=O( X^{1+\epsilon})$, A.~Walfisz
\cite{wa1}: $Q(X)=O( X(\log X)^3)$, and I.~K\'atai \cite{ka}: $
Q(X)=O( X(\log X)^2)$. Later on, E.~Preissmann \cite{Pr} found a
short and elegant proof for this last result, using a deep
inequality of H.L. Montgomery and R.C. Vaughan \cite{mv}.\msk
Recently W.G. Nowak \cite{no4} succeeded in improving K\'{a}tai's
estimate to $$ Q(X) = O\left(X\,(\log X )^{3/2}\,\log_2X\right)\,.
\eqno{(1.15)} $$ We shall sketch the idea of this refinement. On
any interval $[{1\over2} X,X]$, one can split up $P(x)$ as $$
P(x)=H(X,x)+R(X,x)\,, $$ where
$$ \Int_{X/2}^X R(X,x)^2\, {\rm d}x \ll X^{1/2}\,, $$ and
$$ \Int_{X/2}^X H(X,x)^2\, {\rm d}x = C\,\left(X^{3/2}-({\textstyle{1\over2}} X)^{3/2}\right)
+ O\left(X\sum_{1\le n<X^5} {r(n)^2\over n\,\Delta_{{\bf
B}}(n)}\right)\,, \eqno{(1.16)} $$ with
$$ \Delta_{{\bf B}}(n) := \min_{k\in{{\bf B}},\ k\ne n} |k-n|\,,
\qquad {{\bf B}} := \{n\in\Z:\ r(n)>0\ \}\,. $$ This is deduced by
standard techniques, applying (1.6) and a Hilbert-type inequality
of H.L. Montgomery and R.C. Vaughan \cite{mv}. \ssk

Using the trivial bound   $\Delta_{\bf B}(n)\ge1$ in (1.16), along
with the well-known fact that $\sum_{n\le y}r(n)^2\ll y \log y$,
one readily obtains the result of I. K\'atai \cite{ka}, in a way
that mimicks E. Preissmann's \cite{Pr} argument. The basic idea of
the improvement in (1.15) is the observation that the set ${\bf
B}$ is actually rather sparse, i.e., $\Delta_{\bf B}(n)$ should be
"on average" somewhat larger. It was known to E.~Landau \cite{la1}
that
$$ \sum_{1\le n\le x,\ n\in{\bf B}} 1 \ \sim {c\,x\over\sqrt{\log x}}\qquad (x\to\infty)\,, $$
for some $c>0$. Roughly speaking, this suggests that the elements
of ${\bf B}$ behave like some sequence $[c'\,n \sqrt{\log n}]$.
For the latter, the distance between subsequent terms
asymptotically equals $c'\,\sqrt{\log n}$. If one could use this
bound for $\Delta_{\bf B}(n)$ in (1.16), the estimate (1.15) would
be immediate. However, the following result, which can be proved,
and turns out to suffice for this purpose: For every integer
$h\ne0$,
$$ \sum_{1\le n\le x,\ n+h\in{\bf B}} r(n)^2 \ll
\left({|h|\over\phi(|h|)}\right)^2\,x(\log x)^{1/2}\,, \eqno{(1.17)}
$$ where $\phi$ denotes Euler's function. This inequality is
derived by an elementary convolution argument from the related
bound $$ \sum_{1\le n\le x\atop n\in{\bf B},\ kn+h\in{\bf B}} 1 \
\ll \left({|h|k\over\phi(|h|k)}\right)^2\,{x\over\log x} \,,
$$ which is true for arbitrary integers $k>0$ and $h\ne0$. The
latter estimate follows by Selberg's sieve method, and was
established (for $k=1$) by G.J.~Rieger \cite{ri}.\msk We conclude
this subsection by the remark that, during the last 15 years,
higher power moments of $P(x)$ have been investigated as well. We
mention the papers of K.M.~Tsang \cite{Ts1}, D.R.~Heath-Brown
\cite{HB3}, A.~Ivi\'c and P.~Sargos \cite{IS}, and W.~Zhai
\cite{Zh} which provide results of the shape
$$ \Int_0^X (P(x))^m \, {\rm d}x \ \sim \ (-1)^m\,C_m\,X^{1+m/4}\qquad (X\to\infty)\,,  \eqno{(1.18)}  $$
where $3\le m\le9$ and $C_m>0$ are explicitly known constants.
This not only supports the conjecture (1.12), but is of particular
interest in the cases of odd $m$: It means that there is some
excess of the values $x$ for which $P(x)<0$ over those with
$P(x)>0$. This observation matches well with the different
accuracy of the $\Omega_+$- and $\Omega_-$-estimates (1.9) and
(1.10).

\subsubsection*{2.~Lattice points in spheres. }

As most experts agree, in our familiar three-dimensional
Euclidean space the analogue of the Gaussian problem is the most
difficult and enigmatic. For integers $n\ge0$, let
$$ r_3(n) := \#\{(m,k,\ell)\in\Z^3:\ m^2+k^2+\ell^2 = n\ \}\,. $$
Then the formula
$$ \sum_{0\le n\le x} r_3(n) = {4\pi\over3}\,x^{3/2} + P_3(x) \eqno{(2.1)}$$
defines the lattice point number, as well as the lattice point
discrepancy $P_3(x)$, of a compact, origin-centered ball of radius
$\sqrt{x}$. By the same argument which furnished (1.2), it is
trivial to see that $P_3(x)\ll x$. However, it was known already
to E.~Landau \cite{la2} that
$$ P_3(x) = O\left(x^{3/4}\right)\,.  \eqno{(2.2)} $$
For many decades, the subsequent history of improvements of this
$O$-bound was essentially the personal property of
I.M.~Vinogradov, whose series of papers on the topic in 1963
culminated in the result \cite{vi}
$$ P_3(x) = O\left(x^{2/3}(\log x)^6\right)\,.\eqno{(2.3)} $$
Also in this case, the last few years brought significant further
progress in this problem, namely the estimates of F.~Chamizo and
H.~Iwaniec \cite{ci} $$ P_3(x) =
O\left(x^{29/44+\epsilon}\right)\eqno{(2.4)}$$ and of
D.R.~Heath-Brown \cite{hb}
$$ P_3(x) = O\left(x^{21/32+\epsilon}\right)\,.\eqno{(2.5)} $$ Note
that ${29\over44}=0.65909\dots$, ${21\over32}=0.65625$. These
bounds are based on the finer algebro-arithmetic theory of
$r_3(n)$, notably on the explicit formula, an account on which can
be found, e.g., in P.T.~Bateman \cite{ba}. \msk In contrast to
this quite long history of upper bounds, the sharpest known
$\Omega$-result was obtained by G.~Szeg\"o  as early as 1926
\cite{sz} and reads
$$  P_3(x) = \Omega_-\left(x^{1/2}(\log x)^{1/2}\right) \,.\eqno{(2.6)} $$
Neither Hafner's nor Soundararajan's ideas were able to improve
upon this bound, probably because the function $r_3(n)$ is
distributed quite evenly and is not supported on a small subset of
the integers. \msk However, the corresponding $\Omega_+$-estimates
have a longer history of their own: After K.~Chandrasekharan's and
R.~Narasimhan's result \cite{cn} that
$$ \limsup_{x\to\infty} {P_3(x)\over x^{1/2}} = + \infty\,,  $$
W.G. Nowak \cite{no1} proved that
$$ P_3(x) =
\Omega_+\left(x^{1/2}(\log_2x)^{1/2}(\log_3x)^{-1/2}\right)\,.
\eqno{(2.7)} $$ Later on, S.D.~Adhikari and Y.-F.S.~P\'etermann \cite{adh}
obtained the improvement
$$ P_3(x) = \Omega_+\left(x^{1/2}\log_2x\right)\,.\eqno{(2.8)} $$
Finally, and again quite recently, K.M.~Tsang \cite{ts}
complemented Szeg\"o's bound (2.6), showing that in fact
$$  P_3(x) = \Omega_\pm\left(x^{1/2}(\log x)^{1/2}\right) \,.\eqno{(2.9)}
$$ \msk

As for the mean-square of the lattice point discrepancy $P_3(x)$,
its asymptotic behavior has been investigated by V.~Jarnik as
early as in 1940 \cite{ja}. He showed that, with a certain
constant $c_3>0$,
$$ \Int_0^X (P_3(x))^2\, {\rm d}x = c_3\, X^2 \log X + O\left(X^2
(\log X)^{1/2}\right)\,. \eqno{(2.10)}  $$ The proof of (2.10) is
much more difficult than in the planar case. It uses the theory of
theta-functions, complex integration and the classic
Hardy-Littlewood method. As far as the authors were able to
ascertain, an improvement of the error term in (2.10) has never
been attained. \ssk The interesting point of (2.10) is that it
contains the omega estimates (2.6), resp., (2.9), in full
accuracy, apart from the ambiguity of sign. That is, our state of
knowledge is quantitatively quite different than with the circle
problem: While (1.13) together with (1.9) - (1.11) tell us (in the
planar case) that $P(x)\ll x^{1/4}$ {\it in mean-square, with
 unbounded sequences of "exceptional" $x$-values for which
$P(x)$, resp., $-P(x)$, is definitely much larger, } a phenomenon
of this kind is not known for the three-dimensional sphere. \bsk

For spheres of dimensions $s\ge4$, the situation is much easier
and better understood, leaving no room for progress. With
the elementary identity
$$ r_4(n) := \#\{(k_1,k_2,k_3,k_4)\in\Z^4:\ k_1^2+k_2^2+k_3^2+k_4^2=n\ \}
\ =\ 8\,\sigma(n)-32\,\sigma\left({n\over4}\right)\,, $$
$$ \sigma(n)  :=\sum_{m\,|\,n,\ m>0} m   $$ at one's disposal,
it is not difficult to show that (see E.~Kr\"atzel
\cite{kr-anafu}, p.~227-231)
$$ P_4(x):=\sum_{0\le n\le x}r_4(n) - {\pi^2\over2}\,x^2 =
O\left(x\,\log x\right)\,,\quad
P_4(x)=\Omega\left(x\,\log_2x\right)\,. \eqno{(2.11)} $$ These
rather simple estimates have been improved slightly
to
$$P_4(x)= O\left(x\,(\log x)^{2/3}\right)\,,\qquad
P_4(x)=\Omega_\pm\left(x\,\log_2x\right)\,. \eqno{(2.12)} $$ The
$O$-result in (2.12) is due to A.~Walfisz \cite{arnold2},
\cite{wa3} who used some extremely deep and involved analysis,
while the $\Omega$-bounds are due to S.D.~Adhikari and
Y.-F.S.~P\'etermann \cite{adh}. \msk For dimensions $s\ge5$, the
exact order of the lattice point discrepancy $P_s(x)$ of an
origin-centered $s$-dimensional ball of radius $\sqrt x$, has been
known for a long time, namely
$$ P_s(x)= O\left(x^{s/2-1}\right)\,,\qquad
P_s(x)=\Omega\left(x^{s/2-1}\right)\,. \eqno{(2.13)} $$ See again
E.~Kr\"atzel \cite{kr-anafu}, pp 227 ff. For the finer
theory of these higher dimensional lattice point discrepancies,
the reader is referred to the monograph of A.~Walfisz \cite{wa2}.

\input amssym.def
\input amssym
\def\nin{\noindent}
\def\x{{\bf x}}
\def\y{{\bf y}}
\def\u{{\bf u}}
\def\t{{\bf t}}
\def\n{{\bf n}}
\def\0{{\bf 0}}
\def\Re{\mathop{\rm Re}\nolimits}
\def\Im{\mathop{\rm Im}\nolimits}
\def\Res{\mathop{\rm Res}\nolimits}
\def\ind{\mathop{\rm ind}\nolimits}
\newcommand{\upstrut}{{\vrule height 0em depth 0.5em width 0pt}}
\newcommand{\iint}{\int \int}
\newcommand{\idotsint}{\int \cdots \int}
\newtheorem{theorem}{Theorem}
\newtheorem{lemma}{Lemma}
\newtheorem{definition}{Definition}

\subsubsection*{3. Lattice points in convex  bodies.}

By a $p$-dimensional convex domain $CD_p$ we mean a compact convex
set in ${\Bbb{R}}^p$ with the origin as an interior point. Its
distance function $F$ is a homogeneous function of degree $1$,
such that
\[
CD_p = \{ \vec{t} \in {\Bbb{R}^p}: F(\vec{t};CD_p) \le 1 \}.
\]
The number of lattice points $\vec{n} \in {\Bbb{Z}}^p $ { in the
dilated convex domain} $xCD_p$, where $x$ is a large positive
parameter, is denoted by\footnote{Observe the difference in
notation compared to sections 1 and 2: If $CD_p$ is the
$p$-dimensional unit ball, then $A(x;CD_p)=A_p({x^2})$.}
\[
A(x;CD_p) = \# \{ \vec{n} \in {\Bbb{Z}}^p: F(\vec{n};CD_p) \le x \}.
\]
Let $area(CD_2)$ denote the area of $CD_2$ and $vol(CD_p)$ the
volume of $CD_p \; (p > 2)$. The number of lattice points in
$xCD_p$ can be written as
\begin{eqnarray*}
A(x;CD_2) & = & area(CD_2)x^2 + P(x;CD_2),  \\
A(x;CD_p) & = & vol(CD_p)x^p + P(x;CD_p).
\end{eqnarray*}
Here $area(CD_2)x^2$ and $vol(CD_p)x^p$ denote the main terms and
$P(x;CD_p), p \ge 2,$ the remainders or error terms.

\paragraph*{3.1 Lattice points in plane convex domains.} An elementary result due to {M. V. Jarnik}
{ and} { H. Steinhaus} \cite{stein}, 1947: Let $J$ be a closed,
rectifiable { Jordan} curve with area $F$ and length $l
\ge 1$, and let $G$ be the number of lattice points inside and on
the curve. Then $|G - F| < l$.

\smallskip\nin{Conclusion.} Let $l$ denote the length of the curve of boundary
of the convex domain $CD_2$. Then
\[
|P(x;CD_2)| < lx.
\]
The result $P(x;CD_2) \ll x$ is the best of its kind under the
above conditions. An example is the square with length of side $2$
and $(0,0)$ as centre, dilated by $x$.

\medskip\nin{\bf Boundaries with nonzero curvature.} Suppose that
the boundary of the plane convex domain $CD_2$ is sufficiently
smooth with finite nonzero curvature throughout. We say that the
boundary curve is of class $C^k \: (k = 2,3, \ldots )$ if the
radius of curvature is $k - 2$ times continuously differentiable
with respect to the direction of the tangent vector. Let $r$ be
the absolute value of the radius of curvature in a point of the
boundary and $r_{\max}, r_{\min}$ the maximum, minimum of $r$,
respectively. Then we shall assume that
\[
0 < r_{\min} \le r_{\max} < \infty .
\]

\nin{\bf Van der Corput's estimate.} { The theorem of} { J. G. van
der Corput} \cite{corput}, 1920: Let the curve of boundary of the
plane convex domain $CD_2$ be of class $C^2$ , then
\[
P(x;CD_2) \ll x^{\frac{2}{3}}.
\]
{A more precise formulation due to} { H. Chaix} \cite{chaix},
1972:
\[
|P(x;CD_2)| \le C(1 + r_{\max}x)^{\frac{2}{3}}, \quad C > 0\
\hbox{\rm an absolute constant}.
\]
{A formulation with explicit constants is due to} { E. Kr\"atzel}
\cite{kratz5}, 2004:
\[
|P(x;CD_2)| < 48(r_{\max}x)^{\frac{2}{3}} + \left ( 703
\sqrt{r_{\max}} + \frac{3r_{\max}}{5 \sqrt{r_{\min}}} \right )
\sqrt{x} + 11 \quad \mbox{for} \quad x \ge 1.
\]
In particular, if $CD_2 = E_2$ is an ellipse, defined by the
quadratic form
\[
Q(t_1,t_2;E_2) = at_1^2 + 2bt_1t_2 + ct_2^2, \quad ac - b^2 = 1,
\]
then one has
\[
|P(x;E_2)| < 38 x^{\frac{2}{3}} + \left ( 3 a^{\frac{3}{4}} +
700 c^{\frac{3}{4}} \right ) \sqrt{x} + 11.
\]
{E. Kr\"atzel} and { W. G. Nowak} \cite{jewie}, 2004, showed that
the factor $38$ is the leading term of the estimate may be reduced
to $\frac{17}{2}$. In particular, for the circle $E_2 = C \; (a =
c = 1, \: b = 0)$,
\begin{eqnarray*}
|P(x;C)| <
\left\{\begin{array}{ll}
 \frac{17}{2}x^{\frac{2}{3}} + \frac{3}{2} \sqrt{x}
+ 3 & \mbox{for} \quad x \: \ge 28,  \\
9 x^{\frac{2}{3}} &
\mbox{for} \quad x \ge 1000.
\end{array}\right.
\end{eqnarray*}

\nin {V. Jarnik} \cite{jarnik1}, 1925, proved that the estimate of
{ J. G. van der Corput} is the best of its kind. In fact, there
exist plane convex domains $CD_2$ with $P(x;CD_2) = \Omega
(x^{2/3})$.

{J. G. van der Corput} \cite{vandecorp}, 1923, showed that his
estimate can be sharpened for a large class of plane convex
domains $CD_2$. Under the additional assumption that the boundary
curve is of class $C^{\infty}$, he obtained
\[
P(x;CD_2) \ll x^{\vartheta}
\]
holds with some $\vartheta < \frac{2}{3}$.  {O. Trifonov}
\cite{trif}, 1988, showed that the estimate holds with $\vartheta
= \frac{27}{41} + \varepsilon, \varepsilon > 0$.

\medskip\nin{\bf Huxley's estimate.} { The theorem of} { M. N.
Huxley} \cite{Hu1}, 1993: Let the curve of boundary of the plane
convex domain $CD_2$ be of class $C^3$. Then
\[
P(x;CD_2) \ll x^{\frac{46}{73}} (\log x)^{\frac{315}{146}},
\qquad \frac{46}{73} = 0,6301\ldots .
\]
{M. N. Huxley} \cite{Hu2}, 2003, improved this result and
obtained
\[
P(x;CD_2) \ll x^{\frac{131}{208}} (\log x)^{\frac{18637}{8320}},
\qquad \frac{131}{208} = 0,6298\ldots .
\]
{G. Kuba} \cite{kubaell}, 1994, proved a two-dimensional
asymptotic result for a special ellipse $E_2$. Let
\[
A(x,y;E_2) = \# \left \{ (m,n) \in {\Bbb{Z}}^2:
\frac{m^2}{x^2} + \frac{n^2}{y^2} \le 1  \right \} .
\]
Then it is shown that, for $xy \longrightarrow \infty$,
\[
A(x,y;E_2) = \pi xy + P(x,y;E_2)
\]
with

$$
P(x,y;E_2)  \ll
\left\{\begin{array}{ll}
 (xy)^{\frac{23}{73}}(\log
xy)^{\frac{315}{146}} & \mbox{for } x\ll y^{\frac{119}{100}}
(\log xy)^{\frac{315}{100}}, \\
\frac{x}{\sqrt{y}} & \mbox{for } x\gg y^{\frac{119}{100}}(\log
xy)^{\frac{315}{100}}.
\end{array}\right.
$$
In the last case one also has
\[
P(x,y;E_2) = \Omega \left ( \frac{x}{\sqrt{y}} \right )\, .
\]

\bigskip\noindent{\bf Lower bounds.} Suppose that the boundary curve of
the plane convex domain $CD_2$ is of class $C^2$. Then the
following lower bounds have been proved:

\medskip\nin The Theorem of {
V. Jarnik} \cite{jarnik2}, 1924:
\[
P(x;CD_2) = \Omega \left ( x^{\frac{1}{2}} \right ).
\]
{S. Krupi\v{c}ka} \cite{krup}, 1957:
\[
P(x;CD_2) = \Omega_{\pm} \left ( x^{\frac{1}{2}} \right ).
\]
{W. G. Nowak} \cite{now1}, 1985:
\[
P(x;CD_2) = \Omega_{-} \left ( x^{\frac{1}{2}}
(\log x)^{\frac{1}{4}}\right ).
\]
In this connection the {estimates of the average order} of the
remainder are interesting:  \newline { S. Krupi\v{c}ka}
\cite{krup}, 1957:
\[
\int\limits_0^x |P(\sqrt{t};CD_2)|\, {\rm d}t \gg x^{\frac{5}{4}}.
\]
{W. G. Nowak} \cite{now1}, 1985:
\[
\int\limits_0^x |P(\sqrt{t};CD_2)|\, {\rm d}t \ll x^{\frac{5}{4}},
\qquad \int\limits_0^x |P^2(\sqrt{t};CD_2)|\, {\rm d}t \ll
x^{\frac{3}{2}}.
\]
{P. Bleher} \cite{bleher}, 1992:
\[
\int\limits_0^x P^2(t;CD_2) \, {\rm d}t \sim Ax^2.
\]
$A$ is a positive constant depending on $CD_2$.
\par\medskip\nin W. G. Nowak \cite{no3}, 2002:
\[
\int\limits_{x-\Lambda}^{x+\Lambda} P^2(t;CD_2) \, {\rm d}t \sim
4A\Lambda x\, ,
\]
for any $\Lambda=\Lambda(x)$ satisfying
\[\lim_{x\to\infty} {\log x\over \Lambda(x)}=0.\]

\bigskip\nin{\bf Points with curvature zero at the boundary.}

\medskip\nin{\bf Lam\'{e}'s curves.} Let $k \in \Bbb{N}, \: k \ge 3$. A {
Lam\'{e}} curve $L_k$ is defined by
\[
|t_1|^k + |t_2|^k = 1.
\]
The curvature of this curve vanishes in the points $(0,\pm 1), \:
(\pm 1,0)$. The number of lattice points inside and on the dilated
{ Lam\'{e}} curves $xL_k$ is defined by
\[
A(x;L_k) = \# \{ (m,n) \in {\Bbb{Z}}^2: |m|^k + |n|^k \le x^k \} .
\]
The points with curvature zero give an important contribution to the
estimate of the number of lattice points.

\smallskip\nin{\bf A first estimate.} The number of lattice points is
represented by

\[
A(x;L_k) = \frac{2}{k} \frac{\Gamma^2(\frac{1}{k})}
{\Gamma (\frac{2}{k})} \, x^2 + \Delta (x;L_k).
\]
Trivially  $\Delta (x;L_k) \ll x$. The first progress concerning
the upper bound was made by { D. Cauer} \cite{cauer}, 1914, who
proved that
\[
\Delta (x;L_k) \ll x^{1 - \frac{1}{2k-1}},
\]
which was also obtained by { J. G. van der Corput} \cite{holland},
1919. The resolution of the problem of the size of
\[
\Delta (x;L_k) = O, \; \Omega \left ( x^{1 - \frac{1}{k}} \right )
\]
 was given by {B.
Randol} \cite{burton}, 1966, for even $k > 2$ and by {E.
Kr\"atzel} \cite{kr-kreis}, 1967, for odd $k > 3$ and by the same
author \cite{kr-gitter}, 1969, for $k = 3$.

\medskip\nin{\bf The second main term.} {The phenomenon of a second main
term} was introduced by {E. Kr\"atzel} \cite{kr-kreis}, 1967,
\cite{kr-gitter}, 1969 (see also \cite{kr-lp}):
\[
A(x;L_k) = \frac{2}{k} \frac{\Gamma^2(\frac{1}{k})}
{\Gamma (\frac{2}{k})} \: x^2 + 4 \psi_{2/k}^{(k)}(x) + P(x;L_k).
\]
The function $\psi_{\nu}^{(k)}(x)$ is defined for $k \ge 1$, and
for $\nu > \frac{1}{k}$ it is represented by the absolutely
convergent infinite series
$$
\psi_{\nu}^{(k)}(x) = 2 \sqrt{\pi} \, \Gamma \left ( \nu + 1 -
\frac{1}{k} \right ) \sum_{n=1}^{\infty} \left ( \frac{x}{\pi n}
\right )^{\frac{k\nu}{2}} J_{\nu}^{(k)}(2\pi nx),\eqno{(3.1)}
$$
where $J_{\nu}^{(k)}(2\pi nx)$ denotes the generalized {Bessel}
function (see E. Kr\"atzel \cite{kr-lp}, p. 145). It is known that
\[
4 \psi_{2/k}^{(k)}(x) = O, \; \Omega \left ( x^{1 - \frac{1}{k}}
\right ).
\]
Accordingly, this term in the asymptotic representation of
$A(x;L_k)$ is called the {second main term}. The remainder $P$ is
estimated by
\[
P(x;L_k) = o \left ( x^{\frac{2}{3}} \right ) \qquad(x\to\infty)
\]
for $k \ge 3$.

\nin {Improvements:}

\medskip\nin {E. Kr\"atzel} \cite{kr-expo}, 1981 and
independently {W. G. Nowak} \cite{now2}, 1982:
\[
P(x;L_k) \ll x^{\frac{27}{41}}, \qquad \frac{27}{41} = 0,6585 \ldots .
\]
{W. M\"uller} and {W. G. Nowak} \cite{wowe}, 1988:
\[
P(x;L_k) \ll x^{\frac{25}{38}} (\log x)^{\frac{14}{95}}, \qquad
\frac{25}{38} = 0,6578 \ldots .
\]
{G. Kuba} \cite{kuba1}, 1993:
\[
P(x;L_k) \ll x^{\frac{46}{73}}(\log x)^{315\over146}, \qquad \frac{46}{73} = 0,6301 \ldots .
\]

\nin{\bf Lower bounds:}

\medskip\nin {L. Schnabel} \cite{schnabel}, 1982:
\[
P(x;L_k) = \Omega \left ( x^{\frac{1}{3}} \right ).
\]
{E. Kr\"atzel} \cite{kr-lp}, 1988:
\[
P(x;L_k) = \Omega_{\pm} \left ( x^{\frac{1}{2}} \right ).
\]
{W. G. Nowak} \cite{now3}, 1997:
\[
P(x;L_k) = \Omega_- \left ( x^{\frac{1}{2}} ( \log x)^{\frac{1}{4}}
\right ).
\]
{M. K\"uhleitner, W. G. Nowak, J. Schoissengeier} and {T.
Wooley} \cite{kuehleinow1}, 1998:
\[
P(x;L_3) = \Omega_+ \left ( x^{\frac{1}{2}}
( \log \log x)^{\frac{1}{4}} \right ).
\]
{The average order:} {W. G. Nowak} \cite{now4},1996:
\[
\frac{1}{x} \int\limits_0^x P^2(t;L_k) \, {\rm d}t \ll x.
\]
{M. K\"uhleitner} \cite{kuehlei1}, 2000:
\[
\frac{1}{x} \int\limits_0^x P^2(t;L_k) \, {\rm d}t = C_kx + (x^{1
- \varepsilon})
\]
with explicitly given $C_k > 0, \, \varepsilon = \varepsilon (k)
> 0$.

\par\medskip\nin{W. G. Nowak} \cite{now5}, 2000:
\[
\int\limits_{x - 1/2}^{x + 1/2} P^2(t;L_k) \, {\rm d}t \ll x(\log
x)^2.
\]
{M. K\"uhleitner} and {W. G. Nowak} \cite{kuehleinow2}, 2001:
Let $\Lambda (x)$ be an increasing function with
\[
\Lambda (x) \le \frac{x}{2} \quad \mbox{and} \quad
\lim_{x\to\infty} \frac{\log x}{\Lambda (x)} = 0.
\]
Then
\[
\int\limits_{x - \Lambda}^{x + \Lambda} P^2(t;L_k) \, {\rm d}t
\sim 4C_k\Lambda x
\]
with the same constant $C_k$ as above.

\bigskip\nin{\bf General boundary curves.} Suppose that the boundary curve of the plane convex domain $CD_2$ contains a finite
number of points with curvature zero. {Y. Colin de Verdi\`ere}
\cite{colin}, 1977, showed that
\[
P(x;CD_2) \ll x^{1 - \frac{1}{k}}
\]
holds if $k - 2 \ge 1$ is the maximal order of a zero of the
curvature on the boundary of $CD_2$. Moreover, this estimate
cannot be improved if the slope of the boundary curve is rational
in at least one point, where the curvature vanishes to the order
$k - 2$. Further, {B. Randol} \cite{randol}, 1969, {Y. Colin de
Verdiere} \cite{colin}, 1977, and {M. Tarnopolska-Weiss},
\cite{weiss}, 1978, showed: If $CD_2$ is rotated about the origin
by an angle $\varphi$ then for this new domain $CD_2'$ the
estimate
\[
P(x;CD_2') \ll x^{\frac{2}{3}}
\]
holds for almost all $\varphi \in [0,2\pi ]$.

{W. G. Nowak} \cite{now6}, 1984, obtained a refinement in the
rational case. Assuming that the boundary curve has rational slope
in an isolated point with curvature zero, then the contribution of
this point to the discrepancy in the asymptotic representation of
the number of lattice points can be given explicitly by a
{Fourier} series which is absolutely convergent. { E. Kr\"atzel}
\cite{kr-anafu}, 2000, gives an integral representation for this
contribution. It has the precise order $x^{1 - 1/k}$ if $k - 2$ is
the order of vanishing of the curvature in the boundary point. In
other words: Each isolated boundary point with curvature zero and
rational slope produces a new main term in the asymptotic
representation of the number of lattice points.

Further, {W. G. Nowak} \cite{now7}, 1985, obtained a refinement
also in the irrational case also. Assume that the slope in a
boundary point of vanishing curvature is irrational. Under certain
assumptions on the approximability of this irrational number by
rationals, this point gives a contribution to the error term of
order $x^{\vartheta}$ with $\vartheta < \frac{2}{3}$. Improvements
of this result were given by {W. M\"uller} and {W. G. Nowak},
\cite{wolfgeorg}, 1985. Moreover, {M. Peter}, \cite{peter1}, 2000,
proved that there is an error term with $\vartheta < \frac{2}{3}$
if and only if the slope is of finite type at the point with
curvature zero.

\paragraph*{3.2 Lattice points in convex bodies of higher dimensions.} An elementary
result is the Theorem of {J. M. Wills}, \cite{konmenge}, 1973:
Let $K$ be a $p$-dimensional, strictly convex body $(p \ge 3)$ and
$r$ the radius of its greatest sphere in the interior. Let $G$ be
the number of lattice points inside and on the body. Then
\[
vol(K) \left ( 1 - \frac{\sqrt{p}}{2r} \right )^p \le G \le
vol(K) \left ( 1 + \frac{\sqrt{p}}{2r} \right )^p,
\]
provided that $r > \sqrt{p}\, /2$.

\smallskip\nin {Conclusion:}
\[
A(x;CD_p) = vol(CD_p)x^p + O(x^{p - 1}).
\]

\bigskip\nin{\bf Boundaries with nonzero Gaussian curvature.} Suppose that
the boundary of the convex domain $CD_p \: (p \ge 3)$ is a smooth
$(p - 1)$-dimensional surface with finite nonzero {Gauss}ian
curvature throughout. Assume that the canonical map, which sends
every point of the unit sphere in ${\Bbb{R}}^p$ to that point of
the boundary of the convex domain where the outward normal has the
same direction, is one-one and real-analytic. The first
non-trivial estimates for the remainder $P(x;CD_p)$ were given by
{E. Landau}, 1912 and 1924, in case of an ellipsoid and by {E.
Hlawka}, 1950, for a general convex body. We begin with the
general case, and we shall give a short survey of the history of
the ellipsoid problem.

\medskip\nin{\bf The general case.}

\smallskip\nin{The Theorem of} {E. Hlawka}, \cite{hla-gitter},
\cite{hla-omega}, 1950:

\begin{eqnarray*}
P(x;CD_p)=
\left\{\begin{array}{ll}
 O(x^{p - 2 + \frac{2}{p + 1}}),  \\
\Omega( x^{\frac{p - 1}{2}}).
\end{array}\right.
\end{eqnarray*}

\nin{Improvements:}

\medskip\nin{E. Kr\"atzel} and {W. G. Nowak},
\cite{krae-now}, 1991:
\begin{eqnarray*}
P(x;CD_p) \ll
\left\{\begin{array}{ll}
 x^{p - 2 + \frac{12}{7p + 4}} \quad \mbox{for}
\quad 3 \le p \le 7, \\
x^{p - 2 + \frac{5}{3p + 1}} \quad \mbox{for}
\quad 8 \le p.
\end{array}\right.
\end{eqnarray*}
{E. Kr\"atzel} and {W. G. Nowak} \cite{ekk-geo}, 1992:
\begin{eqnarray*}
P(x;CD_p) \ll
\left\{\begin{array}{ll}
 x^{p - 2 + \frac{8}{5p + 2}}( \log x
)^{\frac{10}{5p + 2}} \quad \mbox{for} \quad 3 \le p \le 6, \\
 x^{p - 2 + \frac{3}{2p}}( \log x)^{\frac{2}{p}} \qquad
\quad \mbox{for} \quad 7 \le p.
\end{array}\right.
\end{eqnarray*}
{W. M\"uller} \cite{mueller}, 2000:
\begin{eqnarray*}
P(x;CD_p) \ll
\left\{\begin{array}{ll}
 x^{1 + \frac{20}{43} + \varepsilon} \qquad \qquad
\mbox{for} \quad p = 3, \\
x^{2 + \frac{6}{17} + \varepsilon} \qquad \qquad
\mbox{for} \quad p = 4, \\
x^{p - 2 + \frac{p + 4}{p^2 + p + 2} + \varepsilon} \quad \;
\mbox{for} \quad 5 \le p
\end{array}\right.
\end{eqnarray*}
for every $\varepsilon > 0$.

\medskip\nin{W. G. Nowak} \cite{lattice}, \cite{rest}, \cite{body} showed over
the years 1985 - 1991
\begin{eqnarray*}
P(x;CD_p) =
\left\{\begin{array}{ll}
\Omega_- \left (  x(\log x)^{\frac{1}{3}} \right )
\qquad \qquad \mbox{for} \quad p = 3, \\
\Omega_{\pm} \left ( x^{\frac{p - 1}{2}}(\log x
)^{\frac{1}{2} - \frac{1}{2p}} \right ) \quad \mbox{for} \quad
p \ge 4.
\end{array}\right.
\end{eqnarray*}

\bigskip\nin{A mean square estimate:} {W. M\"uller} \cite{womue}, 1997: If
$p \ge 4$, then
\[
\int\limits_0^x |P(t;CD_p)|^2\, {\rm d}t \ll x^{2p - 3 +
\varepsilon}
\]
for every $\varepsilon > 0$.

\bigskip\nin{\bf Bodies of revolution.}
Let $F(t,z)$ be a distance function depending on two variables
$t, z$. Then
\[
R_p = \left \{ (t_1,t_2, \ldots ,t_{p - 1},z) \in {\Bbb{R}}^p:
F \left ( \sqrt{t_1^2 + t_2^2 + \cdots + t_{p - 1}^2}\, ,z \right )
\le 1 \right \}
\]
denotes a $p$-dimensional body of revolution. {F. Chamizo}
\cite{chamizo}, 1998, proved for the number of lattice points
$A(x;R_p)$ the representation
\[
A(x;R_p) = vol(R_p)x^p + P(x;R_p),
\]
where the remainder is estimated by
\begin{eqnarray*}
P(x;R_p) \ll
\left\{\begin{array}{ll}
  x^{\frac{11}{8}} \: \quad \quad \; \, \mbox{for}
\quad p = 3, \\
 x^3 \log x \quad \mbox{for} \quad p = 5,  \\
 x^{p - 2} \: \qquad \mbox{for} \quad p > 5.
\end{array}\right.
\end{eqnarray*}
{M. K\"uhleitner} \cite{kuehleit}, 2000, proved the  lower bound
\[
P(x;R_3) = \Omega_- \left ( x(\log x)^{\frac{1}{3}} (\log \log
x)^{\frac{\log 2}{3}} e^{-A \sqrt{\log \log \log x}} \right ),
\]
where $A$ is a positive constant.
\par\medskip\nin{M. K{\"u}hleitner and W. G. Nowak} \cite{kn}, 2004, improved
this to
\[
P(x;R_3) = \Omega_- \left ( x(\log x)^{\frac{1}{3}}
(\log \log x)^{\frac{2}{3}(\sqrt{2}-1)}(\log\log\log x)^{-\frac{2}{3}} \right ).
\]

\par\smallskip\nin {An effective estimate in case of} $p = 3$: {E. Kr\"atzel} and {
W. G. Nowak}, \cite{jewie}, 2004, proved a qualitatively weaker
estimate, but with precise numerical constants. We describe only
the simplest case. Let the body of revolution be given by
\[
R_3 = \left \{ (t_1,t_2,z) \in {\Bbb{R}}^3: t_1^2 + t_2^2 \le
f^2(z), \; |z| \le 1 \right \},
\]
where it is assumed that
\[
f(z) \ge 0, \quad f(z) = f(-z), \quad f(1) = 0, \quad f(0) \le 1,
\]
\[
f'(z) \longrightarrow - \infty, \quad f''(z) \longrightarrow
- \infty \quad \mbox{for} \quad z \longrightarrow 1 - 0.
\]
Further assume that $f''(z)$ is monotonic and $0 < r_{\min} \le
r_{\max} < \infty$ for the absolute value of the radius of
curvature. Let
\begin{eqnarray*}
A(x;R_3) & = & \# \left \{ (n_1,n_2,m) \in {\Bbb{Z}}^3:
n_1^2 + n_2^2 \le x^2 f^2 \left( \frac{m}{x} \right ),
\; |m| \le x \right \} \\
& = & vol(R_3)x^3 + P(x;R_3).
\end{eqnarray*}
Then
\[
|P(x;R_3)| \le \left ( 73 r_{\max}^{\frac{3}{4}} +
13 r_{\max}^{\frac{1}{4}} \right ) x^{\frac{3}{2}} + c_1x \log x
+ c_2 x + c_3x^{\frac{3}{4}} + c_4 x^{\frac{1}{2}} + 5
\]
with explicit positive constants $c_1, \ldots , c_4$ depending on
$r_{\min}$ and $r_{\max}$.

\bigskip\nin{\bf Ellipsoids.}
The square of the distance function of an ellipsoid $CD_p = E_p
\; (p \ge 3)$ is given by the positive definite quadratic form
\[
F^2(\vec{t};E_p) = \sum_{\nu = 1}^p \sum_{\mu = 1}^p a_{\nu \mu}
t_{\nu}t_{\mu}, \qquad a_{\nu \mu} = a_{\mu \nu} \in \Bbb{R},
\]
with the determinant $d = det(a_{\nu \mu}) > 0$. Then the
ellipsoid $E_p$ is defined by
\[
E_p = \left \{ \vec{t} \in {\Bbb{R}}^p: F(\vec{t};E_p) \le 1
\right \} ,
\]
and the number of lattice points is given by
\begin{eqnarray*}
A(x;E_p) & = & \# \left \{ \vec{n} \in {\Bbb{Z}}^p:
F(\vec{n};E_p) \le x \right \}  \\
& = & \frac{\pi^{p/2}}{\Gamma (\frac{p}{2} + 1) \sqrt{d}} \, x^p +
P(x;E_p).
\end{eqnarray*}
It is necessary to distinguish between rational and irrational
ellipsoids. An ellipsoid is called rational if there exists a
number $c > 0$ such that $ca_{\nu \mu}$ are integers for all $\nu
, \mu$. An irrational ellipsoid is then a non-rational ellipsoid.

\par\medskip\nin {The Theorem of} {E. Landau} \cite{landau1}, 1915,
\cite{landau2}, 1924: The estimates
\begin{eqnarray*}
P(x;E_p)  =
\left\{\begin{array}{ll}
O\left(x^{p - 2 + \frac{2}{p + 1}}\right), \\
\Omega \left ( x^{\frac{p - 1}{2}} \right )
\end{array}\right.
\end{eqnarray*}
hold for arbitrary ellipsoids.

\medskip\nin{\bf Rational ellipsoids.} {A. Walfisz} \cite{arnold1}, 1924, for $p \ge 8$ and {E. Landau}
\cite{landau3}, 1924, for $p \ge 4$:
\begin{eqnarray*}
P(x;E_p)  \ll
\left\{\begin{array}{ll}
x^{p - 2} \qquad \quad \mbox{for} \quad p > 4, \\
x^2 \log^2x \quad \; \mbox{for} \quad p = 4.
\end{array}\right.
\end{eqnarray*}
Note that
\[
P(x;E_p) = \Omega \left ( x^{p - 2} \right ) \quad \mbox{for} \quad
p \ge 3.
\]
{A. Walfisz} \cite{arnold2}, 1960:
\[
P(x;E_4) = O\left ( x^2(\log x)^{\frac{2}{3}} \right ).
\]
{Y.-K. Lau} and {K. Tsang} \cite{lau}, 2002:
\[
P(x;E_3) = \Omega_{\pm}(x \sqrt{\log x} \, ).
\]

\medskip\nin{\bf Irrational ellipsoids.} In the case of irrational ellipsoids
a lot of special results are available. To obtain them, one
usually imposes special conditions on the coefficients of the quadratic form.
Most of them were
proved by {B. Divi\v{s}, V. Jarnik, B. Novak} and {A. Walfisz}. It
is nearly impossible to list all the different results. For
details, see {F. Fricker's} monograph \cite{fricker}.

\nin The most important results have been obtained for {quadratic forms of
diagonal type}:
\[
F^2(\vec{t};E_p) = \sum_{\nu = 1}^p a_{\nu \nu}t_{\nu}^2 \qquad a_{\nu\nu}>0\ .
\]
{V. Jarnik} \cite{jarnik3}, 1928: The estimate
\begin{eqnarray*}
P(x;E_p) \ll
\left\{\begin{array}{ll}
x^{p - 2} & \mbox{for} \quad p > 4, \\
x^2 \log^2x &  \mbox{for} \quad p = 4
\end{array}\right.
\end{eqnarray*}
holds for arbitrary diagonal quadratic forms.

\par\medskip\nin{V. Jarnik} \cite{jarnik4}, 1929: The estimate
\[
P(x;E_p) = o(x^{p - 2)} \quad \mbox{for} \quad p > 4
\]
holds for irrational, diagonal, quadratic forms.

\par\medskip\nin{V. Jarnik} \cite{jarnik4}, 1929 for $k = 4$ and {A. Walfisz}
\cite{arnold3}, 1927 for $p > 4$: For each monotonic function
$f(x) > 0$ with $f(x) \longrightarrow 0$ for $x \longrightarrow
\infty$ there is an irrational, diagonal, quadratic form such that
\begin{eqnarray*}
P(x;E_p) =
\left\{\begin{array}{ll}
\Omega (x^{p - 2}f(x)) &\mbox{for} \quad p > 4, \\
\Omega (x^2f(x) \log \log x) & \mbox{for} \quad p = 4.
\end{array}\right.
\end{eqnarray*}

\nin{V. Bentkus} and {F. G\"otze} \cite{begoe1}, 1997: The estimate
\[
P(x;E_p) \ll x^{p - 2} \quad \mbox{for} \quad p \ge 9
\]
holds for arbitrary ellipsoids. This result was extended to $p \ge
5$ by {F. G\"otze} \cite{goetze}, 2004.

\nin{V. Bentkus} and {F. G\"otze} \cite{begoe2}, 1999: The estimate
\[
P(x;E_p) = o(x^{p - 2}) \quad \mbox{for} \quad p \ge 9
\]
holds for irrational ellipsoids.

\bigskip\nin{\bf Points with Gaussian curvature zero at the boundary.}

\smallskip\nin{\bf Super Spheres.}
Let $k,p \in \Bbb{N}, \: k \ge 3, \: p \ge 3$. The super spheres
$SS_{k,p}$ are defined by
\[
|t_1|^k + |t_2|^k + \cdots + |t_p|^k = 1.
\]
The {Gauss}ian curvature of these super spheres vanishes for
$t_{\nu} = 0 \: (\nu = 1,2, \ldots ,p)$. Consider the number of
lattice points inside and on the super spheres $xSS_{k,p}$:
\[
A(x;SS_{k,p}) = \# \left \{ \vec{n} \in {\Bbb{Z}}^p:
|n_1|^k + |n_2|^k + \cdots + |n_p|^k \le x^k \right \} .
\]

\nin{\bf A first estimate:} {B. Randol} \cite{burton2}, 1966, for
even $k$ and {E. Kr\"atzel} \cite{kratz6}, 1973, for odd $k$:
\[
A(x;SS_{k,p}) = V_{k,p}x^p + \Delta (x;SS_{k,p}), \qquad
V_{k,p} = \left ( \frac{2}{k} \right )^p
\frac{\Gamma^p (\frac{1}{k})}{\Gamma (1 + \frac{p}{k})},
\]

\begin{eqnarray*}
\Delta (x;SS_{k,p}) =
\left\{\begin{array}{ll}
O\left( x^{p - 2 + \frac{2}{p + 1}}\right) &\mbox{for} \quad k \le p + 1, \\
O, \Omega \left ( x^{(p - 1)(1 - \frac{1}{k})} \right )
&\mbox{for} \quad k > p + 1.
\end{array}\right.
\end{eqnarray*}
Hence the lattice point problem is settled in case when $k > p +
1$, as in the planar case.

\smallskip\nin{\bf An asymptotic representation:} {The Theorem of} {E.
Kr\"atzel}, \cite{kr-lp}, 1988:
\[
A(x;SS_{k,p}) = V_{k,p}x^p + \sum_{r = 1}^{p - 1} H_{k,p,r}(x) +
\Delta_{k,p}(x).
\]
The terms $H_{k,p,r}(x)$ are defined recursively:
\[
\Delta_{k,1}(x) = - 2 \psi (x) = - 2 \left ( x - [x] - \frac{1}{2}
\right ),
\]
\[
H_{k,p,r}(x) = {p \choose r} (p - r)V_{k,p - r} \int\limits_0^x
(x^k - t^k)^{\frac{p - r}{k} - 1} t^{k - 1} \Delta_{k,r}(t) \,
{\rm d}t
\]
for $r = 1,2, \ldots ,p - 1$. Furthermore
\[
H_{k,p,1}(x) = pV_{k,p - 1} \psi_{p/k}^{(k)}(x),
\]
where $\psi_{\nu}^{(k)}(x)$ is defined in (3.1). Estimates:
\[
\psi_{p/k}^{(k)}(x) = O, \, \Omega_{\pm} \left (
x^{(p - 1)(1 - \frac{1}{k})} \right ),
\]
\begin{eqnarray*}
H_{k,p,r}(x) \ll
\left\{\begin{array}{ll}
x^{p - 2} & \mbox{for} \quad p - r > k, \\
x^{(p - r)(1 - \frac{2}{(r + 1)k}) + r - 2 + \frac{2}{r + 1}}
& \mbox{for} \quad p - r \le k,
\end{array}\right.
\end{eqnarray*}
$r = 2,3, \ldots ,p - 1$, and
\[
\Delta_{k,p}(x) \ll x^{p - 2 + \frac{2}{p + 1}}.
\]
\begin{eqnarray*}
H_{k,p,r}(x) & = & \Omega_{\pm} \left (
x^{(p - r)(1 - \frac{1}{k}) + \frac{r - 1}{2}} \right )
\quad\mbox{for} \quad r= 2,3, \ldots ,p - 1, \\
\Delta_{k,p}(x) & = & \Omega_{\pm} \left ( x^{\frac{p - 1}{2}}
\right ).
\end{eqnarray*}

The representation of the number of lattice points suggests that
besides of the main term $V_{k,p}x^p$ there are $p - 1$ further
main terms. But the known upper bounds for these terms at present
allow at most one second main term. So far we have
\begin{eqnarray*}
A(x;SS_{k,p}) =
\left\{\begin{array}{ll}
V_{k,p}x^p + O \left ( x^{p - 2 +
\frac{2}{p + 1}} \right ) &\mbox{for }  k \le p + 1,  \\
V_{k,p}x^p + pV_{k,p - 1} \psi_{p/k}^{(k)}(x) +
O \left ( x^{(p - 2)(1 - \frac{2}{3k}) + \frac{2}{3}} \right ) &\mbox{for } k > p + 1.
\end{array}\right.
\end{eqnarray*}

\nin{Improvements:}

\medskip\nin{E. Kr\"atzel} \cite{kr-lp}, 1988, for $p
= 3$ and {R. Schmidt-R\"oh} \cite{ralph}, 1989 for $p > 3$:
\begin{eqnarray*}
\Delta_{k,p}(x) \ll
\left\{\begin{array}{ll}
x^{p - 2 + \frac{12}{7p + 4}} &
\mbox{for} \quad 3 \le p \le 7, \\
x^{p - 2 + \frac{5}{3p + 1}} &\mbox{for} \quad 8 \le p.
\end{array}\right.
\end{eqnarray*}

\par\medskip\nin{Estimates of the third and fourth main terms:} {E. Kr\"atzel}
\cite{kr-lp}, 1988, for $r = 2$ and {S. H\"oppner} and {E.
Kr\"atzel} \cite{hoeppner}, 1993, for $r = 3$:
\begin{eqnarray*}
\int\limits_0^x \Delta_{k,2}(t^{\frac{1}{k}}) \, {\rm d}t & \ll &
x^{1 - \frac{1}{2k}},  \\
\int\limits_0^x \Delta_{k,3}(t^{\frac{1}{k}}) \, {\rm d}t & \ll &
x \log x
\end{eqnarray*}
and, if $\Delta_{k,r}(t) \ll \Delta_{k,r}^{*}(x)$ for $1 \le t \le x$,
\[
H_{k,p,r}(x) \ll x^{(p - r)(1 + \frac{r - 3}{2k})} \left (
\Delta_{k,r}^{*}(x) \right )^{1 - \frac{p - r}{k}} (\log x)^{r - 2},
\]
provided that $r = 2,3, \: k \ge p - r$. It is highly probable that
the estimate also holds for $r > 3$.

\smallskip\nin{Conclusion:} {G. Kuba's} estimate for $\Delta_{k,2}(x)$ and {
E. Kr\"atzel's} estimate for $\Delta_{k,3}(x)$ lead to
\begin{eqnarray*}
H_{k,p,2}(x) & \ll & x^{(p - 2)(1 - \frac{165}{146k}) + \frac{46}{73}}
(\log x)^{\frac{315}{146}(1 - \frac{p - 2}{k})} \quad
\mbox{for} \quad k \ge p - 2,  \\
H_{k,p,3}(x) & \ll & x^{(p - 3)(1 - \frac{37}{25k}) + \frac{37}{25}}
\qquad \qquad \qquad \qquad \ \mbox{for} \quad k \ge p - 3.
\end{eqnarray*}
In case when $r \ge 4, \: p \ge 5$ it is true that
\[
\sum_{r = 4}^{p - 1} H_{k,p,r}(x) \ll
x^{(p - 4)(1 - \frac{2}{5k}) + \frac{12}{5}}.
\]

\nin{Summary:}

\begin{eqnarray*}
A(x;SS_{k,p}) =
\left\{\begin{array}{ll}
  V_{k,p}x^p + O \left ( x^{\lambda_{k,p}} \right )
& \mbox{for} \quad k < p + 1, \\
 V_{k,p}x^p + pV_{k,p - 1} \psi_{p/k}^{(k)}(x) + O \left (
x^{\vartheta_{k,p}} \right ) & \mbox{for} \quad k \ge p + 1
\end{array}\right.
\end{eqnarray*}
with
\begin{eqnarray*}
\lambda_{k,p} &=&
\left\{\begin{array}{ll}
p - 2 + \frac{12}{7p + 4} & \mbox{for}
\quad 3 \le p \le 7, \\
p - 2 + \frac{5}{3p + 1} & \mbox{for} \quad 8 \le p, \end{array}\right.\\
\vartheta_{k,p} & = & p - 2 + \frac{12}{7p + 4} \quad\ \mbox{for}
\quad p = 3,4, \quad k \le p + 4,
\end{eqnarray*}
\begin{eqnarray*}
\vartheta_{k,p}  =
\left\{\begin{array}{ll}
(p - 2) \left ( 1 - \frac{165}{146k} \right )
+ \frac{46}{73} + \varepsilon & \mbox{for} \quad p = 3,4, \,
k \ge p + 5, \; \varepsilon > 0, \\
 (p - 2) \left ( 1 - \frac{2}{5k} \right ) + \frac{2}{5} +
\frac{4}{5k} & \mbox{for} \quad p \ge 5, \quad \;
k < \frac{533p - 482}{168}, \\
(p - 2) \left ( 1 - \frac{165}{146k} \right ) + \frac{46}{73}
+ \varepsilon & \mbox{for} \quad p \ge 5, \quad \;
k \ge \frac{533p - 482}{168}, \; \varepsilon > 0.
\end{array}\right.
\end{eqnarray*}
In addition to all that, {E. Kr\"atzel} \cite{kratz7}, 1999,
proved that the estimate holds with
\[
\vartheta_{k,p} = p - 2 + \frac{5}{3p + 1}
\]
for even $k$ and $p + 1 \le k < 2p - 4$.

\medskip\nin{\bf Super ellipsoids.}
Let $k,p \in \Bbb{N}, \: k,p \ge 3$. Consider the super ellipsoids
$SE_p$
\[
\lambda_1 |t_1|^k + \lambda_2 |t_2|^k + \cdots + \lambda_p |t_p|^k = 1
\]
with $\lambda_1, \lambda_2, \ldots , \lambda_p > 0$ and the number
of lattice points inside and on the super ellipsoids $xSE_p$:
\[
A(x;SE_p) = \# \left \{ \vec{n} \in {\Bbb{Z}}^p:
\lambda_1 |n_1|^k + \lambda_2 |n_2|^k + \cdots +
\lambda_p |n_p|^k \le x^k \right \}.
\]
{V. Bentkus} and {F. G\"otze} \cite{begoe3}, 2001, proved
\[
A(x;SE_p) = \frac{V_{k,p}}{(\lambda_1 \lambda_2 \cdots
\lambda_p)^{1/k}} \, x^p + P(x;SE_p),
\]
where
\[
P(x;SE_p) \ll x^{p - k}
\]
for even $k$ and sufficiently large $p$.
A super ellipsoid is called rational if there exists a number $c > 0$
such that all $c \lambda_j$ are natural numbers. Otherwise the
superellipsoid is called irrational. Then
\[
P(x;SE_p) = o(x^{p - k})
\]
if and only if $SE_p$ is irrational.

\bigskip\nin{\bf General convex bodies.} Our knowledge about the properties
of points on the boundary with {Gauss}ian curvature zero for the
estimate of the number of lattice points is very poor. Something is
known if the points with {Gauss}ian curvature zero are isolated.

\medskip\nin{K. Haberland} \cite{haber}, 1993: Let the boundary of the convex
body contain only finitely many points of vanishing curvature such
that the tangent plane is rational at these points. Then each such point produces a new main term.

\medskip\nin{E. Kr\"atzel} \cite{nachr}, 2000, and \cite{kr-anafu}, 2000,
simplified the proof of this result for $p = 3$ and gave integral
representations for these contributions. Furthermore, if the slope
of the tangent planes are irrational the error terms will be of
smaller order.

\smallskip The restriction to isolated zeros of curvature excludes some
important bodies such as the super spheres, for example. Thus {M.
Peter}, \cite{peter2}, 2002, extended the considerations to the
case of non-isolated zeros of curvature. Instead of assuming only
finitely many zeros of curvature, he assumed only finitely many
flat points.

The reader is also referred to the papers of {E. Kr\"atzel},
\cite{kratz8}, 2002, and {D.A. Popov}, \cite{popov}, 2000.

\def\DJ{\leavevmode\setbox0=\hbox{D}\kern0pt\rlap
 {\kern.04em\raise.188\ht0\hbox{-}}D}
\font\boldmasi=msbm10 scaled 700      
\def\Bbbi#1{\hbox{\boldmasi #1}}      
\font\boldmas=msbm10                  
\def\Bbb#1{\hbox{\boldmas #1}}        
\def\Qi{{\Bbbi Q}}                      

\subsubsection*{4. Divisor problems and related arithmetic
functions}

The classical (or Dirichlet) divisor problem consists of the
estimation of the function
$$
\Delta(x) := \mathop{\sum\nolimits^{'}}\limits_{n\le x}d(n) -
x(\log x + 2\gamma-1) - {\textstyle{1\over4}}, \eqno{(4.1)}
$$
where $d(n) = \sum_{d|n}1$ is the number of divisors of a natural
number $n$, $\gamma$ is Euler's constant, and the prime
${}^{'}$ denotes that the last summand in (4.1) is halved if $x
\;(\,>1)$ is an integer.

As the reader will observe in the results to come, there is a
far-reaching analogy between this error term $\Delta(x)$ and
$P(x)$, the lattice point discrepancy of the circle, which we
discussed in Section 1. The deep reason is the great similarity of
the generating functions
$$ \sum_{n=1}^\infty {d(n)\over n^s} = \zeta^2(s)
\qquad\hbox{and}\qquad \sum_{n=1}^\infty {r(n)\over n^s} =
4\,\zeta_{\Qi(i)}(s) = 4\, \zeta(s)\,L(s)\quad (\Re(s)>1)\,. $$
Here $\zeta_{\Qi(i)}$ is the Dedekind zeta-function of the
Gaussian field, and $L(s)$ is the Dirichlet $L$-series with the
non-principal character modulo 4.

\smallskip
The {\it generalized (Dirichlet) divisor} problem
(or the Piltz divisor problem, as it is also sometimes called),
consists of the estimation of the function
$$
\Delta_k(x) := \mathop{\sum\nolimits^{'}}\limits_{n\le x}d_k(n) -
xP_{k-1}(\log x) - {\textstyle{(-1)^k\over2^k}}, \eqno{(4.2)}
$$
where $d_k(n)$ is the number of ways $n$ may be written as a
product of $k$ given factors, so that $d_1(n) \equiv 1$ and
$d_2(n) \equiv d(n)$, $\Delta_2(x) \equiv \Delta(x)$. In (4.2),
$P_{k-1}(t)$ is a suitable polynomial in $t$ of degree $k$, and
one has in fact
$$
P_{k-1}(\log x) = \mathop{{\rm
Res}}\limits_{s=1}\,x^{s-1}\zeta^k(s) s^{-1},\eqno{(4.3)}
$$
where $\zeta(s) = \sum_{n\ge1}n^{-s}\;(\Re s>1)$ is
the Riemann zeta-function.
The connection between $d_k(n)$ and $\zeta^k(s)$ is a natural one, since
one has
$$
\zeta^k(s) \;=\; \sum_{n=1}^\infty d_k(n)n^{-s} = \prod_{p\,{\rm
prime}}(1-p^{-s})^{-k}\qquad(\Re s > 1).\eqno{(4.4)}
$$
From (4.4) one infers easily  that $d_k(n)$ is a multiplicative
function of $n$ and that, for a prime $p$ and $\alpha$ a natural
number, $d_k(p^\alpha) = k(k+1)\cdots (k+\alpha-1)/\alpha!.$ The
relation (4.4) may be extended to complex $k$; for this and other
properties of the so-called {\it generalized divisor problem} see
\cite{I13}, Chapter 14. The basic quantities related to
$\Delta_k(x)$ are the numbers
$$
\alpha_k := \inf\,\Bigl\{\;a_k\ge0\;:\;\Delta_k(x)
\ll x^{a_k}\;\Bigr\},\quad
\beta_k := \inf\,\Bigl\{\;b_k\ge0\;:\;\int_1^X\Delta^2_k(x)\,{\rm d}x
 \ll X^{b_k}\;\Bigr\}.\eqno{(4.5)}
$$ Starting from the classical result of Dirichlet that $\Delta(x)
\ll \sqrt{x}$, there have been numerous results on $\alpha_k$ and
$\beta_k$; for some of them the reader is referred to \cite{I13},
Chapter 13, \cite{Ti1}, Chapter 12, \cite{IO}, \cite{kr-lp}. These
results, in principle, have been obtained by two types of methods.
For $k = 2,3$ the estimation of $\alpha_k$ is carried out by means
of exponential sums, and for larger $k$ by employing results
connected with power moments of $\zeta(s)$. The results on
$\beta_k$ have been obtained by several techniques, also using
power moments of $\zeta(s)$. When $k=2$ an important analytic tool
for dealing with $\Delta(x)$ is the formula of G.F. Vorono\"{\i}
(see e.g., \cite{DF}, \cite{I10}, \cite{I13}, \cite{J1},
\cite{V1}). It says that $$ \Delta(x) =
-2\pi^{-1}\sqrt{x}\sum_{n=1}^\infty d(n)n^{-1/2}
\left(K_1(4\pi\sqrt{nx}) +
{\pi\over2}Y_1(4\pi\sqrt{nx})\right),\eqno{(4.6)} $$ where $K_1,
Y_1$ are Bessel functions in standard notation (see e.g.,
\cite{Leb}). Despite the beauty and importance of (4.6), in
practice it is usually expedient to replace it by a truncated
version, obtained by complex integration techniques and the
asympotic formulas for the Bessel functions. This is $$ \Delta(x)
= (\pi\sqrt{2})^{-1}x^{1/4}\sum_{n\le N}d(n)n^{-3/4}
\cos(4\pi\sqrt{nx}-\pi/4) + O_\varepsilon
\left(x^\varepsilon(1+(x/N)^{1/2})\right), \eqno{(4.7)} $$ where
the implied $O$--constant depends only on $\varepsilon$, and the
parameter $N$ satisfies $1\ll N\le x^A$ for any fixed $A>0$.
Estimates of the form
$$ \Delta(x) \ll x^\alpha(\log
x)^C\eqno{(4.8)}
$$
with various values of $\alpha$ and $C\ge0$ have appeared over the
last hundred years or so and, in general, reflect the progress of
analytic number theory. The last in a long line of records (see
\cite{I13}, Chapter 13 for a discussion) for bounds of the type (4.8)
 is $\alpha = {131\over416} = 0.3149\ldots\,,$ due to M.N. Huxley
\cite{Hu2}. This result is obtained by the intricate use of
exponential sum techniques (see his monograph \cite{Hu1})
connected to the Bombieri--Iwaniec method of the estimation of
exponential sums (see the works of E. Bombieri -- H. Iwaniec
\cite{BI} and of H. Iwaniec -- C.J. Mozzochi \cite{IM}). The limit
of these methods appears to be $\alpha_2 \le 5/16$, whilst
traditionally one conjectures that $\alpha_2 = 1/4$ holds. In
general, one conjectures that $\alpha_k = (k-1)/(2k)$ and $\beta_k
= (k-1)/(2k)$ holds for every $k\ge2$. Either of these statements
(see \cite{I13} or \cite{Ti1}) is equivalent to the Lindel\"of
hypothesis that $\zeta({1\over2} + it) \ll_\varepsilon
|t|^\varepsilon$ for any given $\varepsilon>0$; note that
trivially $\beta_k \le \alpha_k$ holds for $k\ge2$.

The best known upper bound $\alpha_3 \le 43/96$ was obtained by G.
Kolesnik \cite{Kol1}, who used a truncated formula for
$\Delta_3(x)$, analogous to (4.7), coupled with the estimation of
relevant three-dimensional exponential sums. Further sharpest
known bounds for $\alpha_k$, obtained by using power moments for
$\zeta(s)$ are (see \cite{I13}, Chapter 13): $\alpha_k \le
(3k-4)/(4k)\;(4 \le k \le 8)$, $\alpha_9 \le 35/54$ and (see
\cite{IO} for the values of $k$ between 10 and 20), $\alpha_k \le
(63k-258)/(64k)\;(79 \le k \le 119)$ (see \cite{IO}). For large
values of $k$ one has the best bounds which come from the
Vinogradov--Korobov zero-free region for $\zeta(s)$ and the bound
for $\zeta(1+it)$. Namely if $D>0$ is such a constant for which
one has
$$
\zeta(\sigma + it) \ll |t|^{D(1-\sigma)^{3/2}}\log^{2/3}|t|\qquad(
{\textstyle{1\over2}} < \sigma \le 1),\eqno{(4.9)}
$$
then
$$
\alpha_k \;\le\;1 - {\textstyle{1\over3}}\cdot 2^{2/3}(Dk)^{-2/3}.
\eqno{(4.10)}
$$
From the work of H.-E. Richert \cite{Ri2} it is known that $D\le
100$ may be taken in (4.9). The best known result $D \le 4.45 $ is
due to K. Ford \cite{Fo}.

\medskip
As for results on $\beta_k$, it is worth noting that the exact
value $\beta_k = (k-1)/(2k)$ is known for $k = 2,3$ (see
\cite{I13}, Chapter 13) and $k=4$ (see \cite{HB1}). One has
$\beta_5 \le 9/20$ (see W. Zhang \cite{Zn}), $\beta_6 \le 1/2,
\beta_7 < 0.55469, \beta_8 < 0.606167, \beta_9 < 0.63809,
\beta_{10} < 0.66717$ (see \cite{IO}). If (4.9) holds, then
(ibid.)
$$
\beta_k \;\le\;1 - {\textstyle{2\over3}}(Dk)^{-2/3}. \eqno{(4.11)}
$$

\medskip
On the other hand, one may ask for lower bounds or omega results
for $\Delta_k(x)$. G.H. Hardy showed in \cite{Ha1} that
$$
\Delta(x) = \cases{\Omega_+((x\log x)^{1/4}\log_2x),&\cr
\Omega_-(x^{1/4}).\cr}
$$
Further improvements are due to K.S. Gangadharan \cite{Ga}, K.
Corr\'adi-- I. Katai \cite{CK} and to J.L. Hafner \cite{Hf1}. K.
Corr\'adi-- I. Katai proved
$$
\Delta(x) =
\Omega_-\left(x^{1/4}\exp\{c(\log_2x)^{1/4}(\log_3x)^{-3/4}\}
\right)\qquad(c>0)\, .
$$
Hafner op. cit. showed that, with suitable constants $A>0, B>0$,

$$
\Delta(x) = \Omega_+\left((x\log x)^{1/4}
(\log_2x)^{(3+\log4)/4}\exp(-A\sqrt{\log_3x})\right)\, .
\eqno{(4.12)}
$$
Recently K. Soundararajan \cite{So} (see his lemma in Section 1.2)
proved
$$
\Delta(x) = \Omega\left((x\log x)^{1/4}
(\log_2x)^{(3/4)(2^{4/3}-1)}(\log_3x)^{-5/8}\right)\, .\eqno{(4.13)}
$$

Note that $(3/4)(2^{4/3}-1) = 1.1398\ldots...\,$, whilst
$(3+\log4)/4 = 1.0695\ldots...\,$, but on the other hand (4.13) is
an omega result and not an $\Omega_+$ or $\Omega_-$-result. In
other words, Soundararajan obtains a better power of $\log_2x$,
but cannot ascertain whether $\Delta(x)$ takes large positive, or
large negative values. His method, like Hafner's (see \cite{Hf2}),
carries over to $\Delta_k(x)$. Hafner proved, with suitable
$A_k>0$,
$$
\Delta_k(x) = \Omega^*\left((x\log x)^{(k-1)/(2k)}
(\log_2x)^{((k-1)/2k)(k\log k -k+1)+k-1}\exp(-A_k\sqrt{\log_3x})\right),
$$
where $\Omega^* = \Omega_+$ if $k=2,3$ and $\Omega^* = \Omega_\pm$ if
$k\ge4$. Soundararajan showed that
$$
\Delta_k(x) = \Omega\left((x\log x)^{(k-1)/(2k)}
(\log_2x)^{((k+1)/2k)(k^{2k/(k+1)}-1)}(\log_3x)^{-1/2-(k-1)/4k}\right).
$$
The above estimate holds with $\Omega_+$ in place of $\Omega$ if
$k\,\equiv\,3\,({\rm  mod}\,8)$, and with $\Omega_-$ in place of
$\Omega$ if $k\,\equiv\,7\,({\rm  mod}\,8)$.

\medskip
Large values and power moments of $\Delta(x)$ were investigated by
A. Ivi\'c  \cite{I3}. It was shown there that
$$
\int_0^X|\Delta(x)|^A\,{\rm d}x \;\ll_\varepsilon\;
X^{1+A/4+\varepsilon}\qquad(0 \le A \le 35/4),\eqno{(4.14)}
$$
and the range for $A$ for which the above bounds holds can be
slightly extended by using newer bounds for $\zeta({1\over2}+it)$.
The bounds in (4.14) are precisely what one expects to get if the
conjectural bound $\alpha_2 = 1/4$ holds. They were used by D.R.
Heath-Brown \cite{HB3} to prove that the function
$x^{-1/4}\Delta(x)$ has a distribution function and that, for
$A\in [0,9]$ (not necessarily an integer), the mean value
$$
X^{-1-A/4}\int_0^X|\Delta(x)|^A\,{\rm d}x
$$
converges to a finite limit as $X\to\infty$. Moreover the same is true
for the odd moments
$$
X^{-1-A/4}\int_0^X\Delta(x)^A\,{\rm d}x \qquad(A = 1,3,5,7,9).
$$
For particular cases sharper results on moments are, in fact,
known. A classical result of G.F. Vorono\"{i} \cite{V1} states
that
$$
\int_1^X\Delta(x)\,{\rm d}x \;=\; {1\over4}X + O(X^{3/4}),
$$
so that $\Delta(x)$ has 1/4 as mean value.
K.C. Tong \cite{To1} proved the mean square formula
$$
\int_1^X\Delta^2(x)\,{\rm d}x \;=\; CX^{3/2} + F(X)
\quad(C = (6\pi^2)^{-1}\sum_{n=1}^\infty d^2(n)n^{-3/2} = 0.6542869\ldots\,),
$$
where the error term $F(X)$ satisfies $F(X) \ll X\log^5X$. Much
later E. Preissmann \cite{Pr} improved this to $F(X) \ll
X\log^4X$. Further relevant results on $F(X)$ are to be found in
the works of Y.-K. Lau- K.M. Tsang \cite{LT} and K.-M. Tsang
\cite{Ts1}, \cite{Ts2}. In particular, it was shown that
$$
\int_2^X F(x)\,{\rm d}x \;=\; -(8\pi^2)^{-1}X^2\log^2X + cX^2\log
X + O(X)\eqno{(4.15)}
$$
holds with a suitable $c>0$. On the basis of (4.15) it is
plausible to conjecture that, with a certain $\kappa$, one has
$$
F(x) = -(4\pi^2)^{-1}x\log^2x + \kappa x\log x + O(x).
$$
K.-M. Tsang \cite{Ts3} treated the third and fourth moment
of $\Delta(x)$, proving
$$
\int_1^X\Delta^3(x)\,{\rm d} x = BX^{7/4} +
O_\varepsilon(X^{\beta+\varepsilon})  \qquad(B > 0)\eqno{(4.16)}
$$
and
$$
\int_1^X\Delta^4(x)\,{\rm d} x = CX^2 +
O_\varepsilon(X^{\gamma+\varepsilon})  \qquad(C > 0)\eqno{(4.17)}
$$
with $\beta = {47\over28}, \gamma  = {45\over23}$. In a
forthcoming work of  A. Ivi\'c and P. Sargos \cite{IS}, those
values are improved to $\beta = {7\over5}, \gamma = {23\over12}$.
A result on integrals of $\Delta^3(x)$ and $\Delta^4(x)$ in short
intervals, which improves a result of W.G. Nowak \cite{No1}, is
also obtained in \cite{IS}. The analogues of (4.16) and (4.17)
hold for the moments of $P(x)$, with the same values of the
exponents $\beta$ and $\gamma$. One of the ingredients used
therein for the proof of (4.17) is the following lemma, due to O.
Robert-- P. Sargos \cite{RS}: {\it Let $k\ge 2$ be a fixed integer
and $\delta > 0$ be given. Then the number of integers
$n_1,n_2,n_3,n_4$ such that $N < n_1,n_2,n_3,n_4 \le 2N$ and}
$$
|n_1^{1/k} + n_2^{1/k} - n_3^{1/k} - n_4^{1/k}| < \delta N^{1/k}
$$
{\it is, for any given $\varepsilon>0$,}
$$
\ll_\varepsilon N^\varepsilon(N^4\delta + N^2).\eqno{(4.18)}
$$
Moments of $\Delta^k(x)$ when $5\le k\le 9$ are treated by W. Zhai \cite{Zh},
who obtained asymptotic formulas with error terms.

\smallskip
Finally we mention results for the Laplace transform of $\Delta(x)$
and $P(x)$, studied in \cite{I8} and \cite{I12}. We have
$$
\int_0^\infty P^2(x)e^{-x/T}\,{\rm d} x
={1\over4}\left({T\over\pi}\right) ^{3/2}\sum_{n=1}^\infty
r^2(n)n^{-3/2} - T +
O_\varepsilon(T^{\alpha+\varepsilon})\eqno{(4.19)}
$$
and
$$\int_0^\infty \Delta^2(x){\rm e}^{-x/T}\,{\rm d} x =
 {1\over8}\left({T\over2\pi}\right)^{3/2}
\sum_{n=1}^\infty d^2(n) n^{-3/2} + T(A_1\log^2T + A_2\log T +
A_3) + O_\varepsilon(T^{\beta+\varepsilon}). \eqno{(4.20)}
$$

The $A_j$'s are suitable constants ($A_1 =
-1/(4\pi^2)$), and the constants ${1\over2} \le \alpha < 1$
and ${1\over2} \le \beta < 1$ are defined by the asymptotic formula
$$
\sum_{n\le x}r(n)r(n+h) = {(-1)^h8x\over h}\sum_{d|h}(-1)^dd +
E(x,h),\; E(x,h) \ll_\varepsilon
x^{\alpha+\varepsilon},\eqno{(4.21)}
$$
$$
\sum_{n\le x}d(n)d(n+h) = x\sum_{i=0}^2(\log
x)^i\sum_{j=0}^2c_{ij} \sum_{d|h}\left({\log d\over d}\right)^j +
D(x,h), \; D(x,h) \ll_\varepsilon
x^{\beta+\varepsilon}.\eqno{(4.22)}
$$
The $c_{ij}$'s are certain absolute constants, and the
$\ll$--bounds both in (4.21) and in (4.22) should hold uniformly
in $h$ for $1\le h \le x^{1/2}$. With the values $\alpha = 5/6$ of
D. Ismoilov \cite{Ism} and $\beta = 2/3$ of Y. Motohashi \cite{Mo}
it followed then that (4.19) and (4.20) hold with $\alpha = 5/6$
and $\beta = 2/3$. Motohashi's fundamental paper (op. cit.)  used
the powerful methods of spectral theory of the non-Euclidean
Laplacian. A variant of this approach was used recently by T.
Meurman \cite{Me} to sharpen Motohashi's bound for $D(x,h)$ for
`large' $h$, specifically for $x^{7/6} \le h \le
x^{2-\varepsilon}$, but the limit of both methods is $\beta = 2/3$
in (4.22). Using the results of F. Chamizo \cite{Ch}, A. Ivi\'c \cite{I12} obtained $\alpha = 2/3$ in (4.19), which appears
to be the limit of the present methods.

\bigskip
The {\it general divisor problem} can be defined in various ways.
Here we shall follow the notation introduced in \cite{I5}. Let
$d(a_1,a_2,\ldots a_k;n)$ be the number of representations of an
integer $n\ge1$ in the form $n = m_1^{a_1}m_2^{a_2}\ldots
m_k^{a_k}$, where  $k\ge2$  and $1 \le a_1 \le a_2 \le \cdots \le
a_k$ are given integers, and the $m$'s are positive integers. Then
by the general divisor problem we shall mean the estimation of the
quantity

\renewcommand{\theequation}{4.23}
\begin{eqnarray}
\Delta (a_1,a_2,\ldots a_k;n) &=& \sum_{n\le x}d(a_1,a_2,\ldots
a_k;n) - \sum_{j=1}^k \,\mathop{{\rm
Res}}\limits_{s=1/a_j}\,\left(\prod_{r=1}^k
\zeta(a_rs)\right)x^ss^{-1}\nonumber \\ &=& \sum_{n\le
x}d(a_1,a_2,\ldots a_k;n) - \sum_{j=1}^k \left(\prod_{r=1,r\ne
j}^k\zeta(a_r/a_j)\right)x^{1/a_j}
\end{eqnarray}

\nin if the $a_j$'s are distinct. If this is not so, then the
appropriate limit has to be taken in the above sum. For instance,
if
 $ a_1= a_2 = \ldots = a_k =1$, then
we get the classical Dirichlet divisor problem (without the dash
in the sum and the constant term in (4.2)). It should be noted
that $k$ does not have to be finite. A good example for this case
is
$$
\zeta(s)\zeta(2s)\zeta(3s)\ldots = \sum_{n=1}^\infty a(n)n^{-s}
\qquad(\Re s > 1),\eqno{(4.24)}
$$
corresponding to $k = \infty,\, a_j = j$ for every $j$. Using the
product representation (4.3) (with $k=1$) for $\zeta(s)$, it
follows that $a(n)$ is a multiplicative function of $n$, and that
for every prime $p$ and every natural number $\alpha$ one has
$a(p^\alpha) =P(\alpha)$, where $P(\alpha)$ is the number of
(unrestricted) partitions of $\alpha$. Thus $a(n)$ denotes the
number of nonisomorphic Abelian (commutative) groups with $n$
elements (see \cite{Krae3} for a survey of results up to 1982).

Another example of this kind is
$$
\prod_{r\ge1,m\ge1}\zeta(rm^2s) = \sum_{n=1}^\infty
S(n)n^{-s}\qquad(\Re s >1),
$$
where $S(n)$ denotes the number of nonisomorphic semisimple rings with
$n$ elements (see J. Knopfmacher \cite{Kno}).

In the general divisor problem we introduce two constants:

\begin{eqnarray*}
{\bar\alpha}_k := \alpha(a_1,a_2,\ldots,a_k) &=&
\inf\left\{\,\alpha\ge0 : \Delta(a_1,a_2,\ldots,a_k) \ll x^\alpha
\right\},\\ {\bar\beta}_k := \beta(a_1,a_2,\ldots,a_k) &=&
\inf\left\{\,\beta\ge0 : \int_1^X\Delta^2(a_1,a_2,\ldots,a_k){\rm
d}x \ll X^{1+2\beta}\right\},
\end{eqnarray*}
which generalize $\alpha_k,\,\beta_k$ in (4.5). From the classical
results of E. Landau \cite{La} it follows that
$$
{k-1\over2(a_1+a_2+\ldots+a_k)} \le {\bar \alpha}_k \le
{k-1\over(k+1)a_1} \eqno{(4.25)}
$$
if the numbers $a_j$ are distinct. In many particular cases the
bounds in (4.25) may be superseded by the use of various
exponential sum techniques, coupled with complex integration
methods and the results  on $\zeta(s)$ (moments, functional
equations etc.). The case $k = 3$ (the three-dimensional divisor
problem) is extensively discussed in  E. Kr\"atzel's monograph
\cite{kr-lp} and e.g., the work of H.-Q. Liu \cite{Li2}.

\smallskip
In what concerns results on ${\bar\beta}_k$, A. Ivi\'c \cite{I5}
proved the following result: let $r$ be the largest integer which
satisfies $(r-2)a_r \le a_1 + a_2 + \ldots + a_{r-1}$ for $2\le r
\le k$, and let
$$
g_k = g(a_1,a_2,\ldots,a_k) := {r-1\over2(a_1+a_2+\ldots+a_r)}.
\eqno{(4.26)}
$$
Then ${\bar\beta}_k\ge g_k$, and if
$$
\int_0^T|\zeta({\textstyle{1\over2}}+it)|^{2k-2}\,{\rm d}t
\ll_\varepsilon T^{1+\varepsilon}\eqno{(4.27)}
$$
holds, then ${\bar\beta}_k = g_k$. Thus assuming (4.27) we obtain
a precise evaluation of ${\bar\beta}_k$, but it should be remarked
that (4.27) at present is known to hold for $k=2,3$, whilst its
truth for every $k$ is equivalent to the Lindel\"of hypothesis.
Moreover, we have
$$
\int_1^X\Delta^2(a_1,a_2,\ldots,a_k)\,{\rm d}x  = \Omega
\left(X^{1+2g_k}\log^AX\right),\eqno{(4.28)}
$$
where $A = A(a_1,a_2,\ldots,a_k)$ is explicitly evaluated in \cite{I5}.
The results remain valid if $k = \infty$, provided that the generating
function is of the form
$$
\zeta^{q_1}(b_1s)\zeta^{q_2}(b_2s)\zeta^{q_3}(b_3s)\ldots
\quad(1\le b_1 < b_2 < b_3<\ldots),\eqno{(4.29)}
$$
where the $b$'s and the $q$'s are given natural numbers. The
generating function (4.24) of $a(n)$ is clearly of the form
(4.29). In this case it is convenient to define the error term as
$$
E(x) := \sum_{n\le x}a(n) - \sum_{j=1}^6c_jx^{1/j} \quad\left(c_j
= \prod_{\ell=1,\ell\ne j}^\infty \zeta\bigl({\ell\over
j}\bigr)\right). \eqno{(4.30)}
$$
Then (4.28) implies that
$$
\int_1^XE^2(x) \,{\rm d}x  = \Omega(X^{4/3}\log X),\eqno{(4.31)}
$$
and (4.31)  yields further $E(x) = \Omega(x^{1/6}\log^{1/2}x)$,
obtained by R. Balasubramanian--K. Ramacha-\break ndra \cite{BR}.
The omega result (4.31) is well in tune with the upper bound
$$
\int_1^XE^2(x) \,{\rm d}x  = O(X^{4/3}\log^{89}X)
$$
of D.R. Heath-Brown \cite{HB2}, who improved the exponent $39/29
= 1.344827\ldots\,$ of \cite{I4}.

\smallskip
As for upper bounds for $E(x)$ of the form $E(x) =
O(x^c\log^Cx)\;(C\ge0)$ or $E(x) =
O_\varepsilon(x^{c+\varepsilon})$, they have a long and rich
history (note that $c\ge1/6$ must hold by (4.31)). The bound
$c\le1/2$ was obtained by P. Erd\H os--G. Szekeres \cite{ES}, who
were the first to consider the function $a(n)$ and the so-called
{\it powerful numbers} (see \cite{I3}, \cite{kr-lp} for an
account). After their work, the value of $c$ was decreased many
times (chronologically in the works \cite{KR}, \cite{Ri1},
\cite{Scw}, \cite{Sch}, \cite{Sr}, \cite{Kol2} and Liu \cite{Li1},
who had $c\le 50/199 = 0.25125\ldots\,)$. Recently O. Robert--P.
Sargos \cite{RS} obtained $c \le 1/4+\varepsilon$  by the use of
(4.18). This is the limit of the method, and the result comes
quite close to the conjecture of H.-E. Richert \cite{Ri1}, made in
1952, that $E(x) = o(x^{1/4})$ as $x\to\infty$. The most
optimistic conjecture is that $c = 1/6+\varepsilon$ holds.

\medskip
A related problem is the estimation of $T(x) = \sum\tau(G)$, where
$\tau(G)$ denotes the number of direct factors of a finite Abelian group
$G$, and summation is over all Abelian groups whose orders do not
exceed $x$. It is known (see \cite{Co} or \cite{Krae4}) that
$$
T(x) = \sum_{n\le x}t(n),\quad\sum_{n=1}^\infty t(n)n^{-s}
= \zeta^2(s)\zeta^2(2s)\zeta^2(3s)\ldots\;(\Re s>1),
$$
which is of the form (4.29) with $q_j = 2, b_j = j$. Here one
should define appropriately the error term (note that $t(n)
= \sum_{d|n} a(d)a(n/d)$) as

\begin{eqnarray*}
\Delta_1(x) &:=& T(x) - \sum_{j=1}^5(D_j\log x + E_j)x^{1/j}\\& =&
\sum_{mn\le x}a(m)a(n)- \sum_{j=1}^5(D_j\log x + E_j)x^{1/j},
\end{eqnarray*}

\nin where $D_j$ and $E_j\,(D_1>0)$ can be explicitly evaluated
(see H. Menzer--R. Seibold \cite{MS}).  They proved that
$\Delta_1(x) = O_\varepsilon(x^{\rho+\varepsilon})$ with $\rho =
45/109 = 0.412844\ldots\,$, improving a result of E. Kr\"atzel
\cite{Krae4}, who had the exponent 5/12. As for the true order of
$\Delta_1(x)$, one expects $\rho = 1/4$ to hold. This is supported
by the estimate (see \cite{I6} and \cite{I7})
$$
\int_1^X\Delta^2_1(x)\,{\rm d}x = \Omega(X^{3/2}\log^4X)\, ,
\eqno{(4.32)}
$$
which yields $\rho \ge 1/4$.
The value of $\rho$ was later improved (see \cite{HM} and \cite{Li2}),
and currently the best bound  $\rho \le 47/130$ is due to J. Wu \cite{Wu}.
Moreover, in \cite{I6} it was shown that
$$
\int_1^X\Delta^2_1(x)\,{\rm d}x = O_\varepsilon
(X^{8/5+\varepsilon}),
$$
while if (4.27) holds with $k=5$, then one has
$$
\int_1^X\Delta^2_1(x)\,{\rm d}x =
O_\varepsilon(X^{3/2+\varepsilon}).\eqno{(4.33)}
$$
In this case (4.32) and (4.33) determine fairly closely the true
order of the mean square of $\Delta_1(x)$.

\medskip
The previous discussion centered on {\it global} problems
involving arithmetic functions which arise in connection with
divisor problems, that is, the estimation of the error terms in
the asymptotic formulas for summatory functions or the related
power moment estimates. One can, of course, treat also {\it local}
problems as well, namely problems involving  pointwise estimation,
distribution of values, and other arithmetic properties. The
literature on this subject, which can be also considered as a part
of the theory of divisor problems, is indeed vast, especially on
$d(n)$ and $d_k(n)$. Thus it would be a great increase in the
length of this work, as well somewhat outside the mainstream of
the theory, if we dwelt in general on local problems. Also one
sometimes, under divisor problems, considers divisor problems of
the type (4.29) (or (4.23)), where the $q$'s are allowed to be
negative. This greatly increases  the class of functions
that are allowed (like e.g., {\it squarefree numbers} whose
characteristic function is generated by $\zeta(s)/\zeta(2s)$, or
{\it squarefull numbers} whose characteristic function is
generated by $\zeta(2s)\zeta(3s)/\zeta(6s)$). In this work we have
found it best to adhere to divisor problems of the form (4.23) or
(4.29).

\medskip
As for local problems, we shall conclude with just a few words on
a representative divisor problem, namely the function $a(n)$. It
was proved by E. Kr\"atzel \cite{Krae1} that
$$
\limsup_{n\to\infty}\,\log a(n)\cdot{\log\log n\over\log n} =
{\log5\over4},\eqno{(4.34)}
$$
and his result was sharpened by W. Schwarz -- E. Wirsing \cite{SW}
and (on the Riemann Hypothesis) by J.-L. Nicolas \cite{NI}.
General results for multiplicative functions, analogous to (4.34),
were obtained by E. Heppner \cite{He} and A.A. Drozdova -- G.A.
Freiman (see \cite{Po}, Chapter 12). The existence of the local
densities of $a(n)$, the numbers $$ d_k :=
\lim_{x\to\infty}\,x^{-1}\sum_{n\le x, a(n)=k}1
$$
for any given integer $k\ge1$, was proved by Kendall-Rankin
\cite{KR}. This was sharpened in \cite{I1} to
$$
\sum_{n\le x, a(n)=k}1 = d_kx + O(x^{1/2}\log x),
$$
where the $O$--constant is uniform in $k$. Further sharpenings and
generalizations are to be found in \cite{I2}, \cite{IT},
\cite{Krae2}, \cite{KW} and \cite{no2}.

\medskip
The distribution of values of $a(n)$ was investigated in \cite{I2}
and P. Erd\H os -- A. Ivi\'c \cite{EI1}. It was shown there that
$[x,\,2x]$ contains at least $\sqrt{x}$ integers $n$ for which
$a(n+1) = a(n+2) = \ldots = a(n+k)$ with $k = [\log
x\log_3x/(40(\log_2x)^2)]$, and at least $\sqrt{x}$ integers $m$
such that the values $a(m+1), a(m+2), \ldots, a(m+t)$ are all
distinct, where for a suitable $C>0$
$$
t = [C(\log x/\log\log x)^{1/2}].
$$
Iterates of $a(n)$ were treated by P. Erd\H os -- A. Ivi\'c
\cite{EI2} and \cite{I7}.

\par\bigskip

\vbox{Aleksandar Ivi\'c \smallskip

Katedra Matematike RGF-a

Universiteta u Beogradu

Djusina 7

11000 Beograd

Serbia, Yugoslavia \smallskip

E-mail: {\tt \ ivic@rgf.bg.ac.yu}

\bigskip\bigskip

Ekkehard Kr\"atzel \smallskip

Institut f\"ur Mathematik

Universit\"at Wien

Nordbergstra\ss e 15

1090 Vienna, Austria \smallskip

http://www.mat.univie.ac.at/\~{}baxa/kraetzel.html

\bigskip\bigskip

Manfred K\"uhleitner and Werner Georg Nowak \smallskip

Institut f\"ur Mathematik

Department f\"ur Integrative Biologie

Universit\"at f\"ur Bodenkultur Wien

Gregor Mendel-Stra{\ss}e 33

1180 Vienna, Austria \smallskip

E-mail: {\tt \ kleitner@edv1.boku.ac.at,\quad nowak@mail.boku.ac.at} \smallskip

http://www.boku.ac.at/math/nth.html}

\end{document}